\documentclass{amsart}
\usepackage{epsf}
\usepackage{graphicx}
\newcommand {\cal} {\mathcal}

\newcommand {\Reals} {\bf R}
\newcommand {\Rtwo} {\Reals^2}
\newcommand {\Rthree} {\Reals^3}
\newcommand {\D} [1] {\nabla_{#1}}
\newcommand {\DD} [1] {{\overline {\nabla}}_{#1}}
\newcommand {\DT} [1] {\tilde D_{#1}}

\newcommand{\tr}{{\rm tr}}

\newcommand {\thickness} {h_0}

\newcommand{\ggg}{\mbox{${\rm g}$}}

\newcommand {\Pe} {\Phi_\epsilon}
\newcommand {\pe} {\phi_\epsilon}
\newcommand {\Ne} {n_\epsilon}

\newcommand {\Tom} {T^{\left( \omega \right)}}
\newcommand {\domega} [1]  {\partial_{#1} \omega}
\newcommand {\dpe} [1]  {\partial_{#1} \pe}
\newcommand {\dPe} [1]  {\partial_{#1} \Pe}
\newcommand {\dphi} [1]  {\partial_{#1} \phi}
\newcommand {\dpsi} [1]  {\partial_{#1} \psi}
\newcommand {\dPhi} [1]  {\partial_{#1} \Phi}
\newcommand {\dn} [1]  {\partial_{#1} n}
\newcommand {\dne} [1]  {\partial_{#1} \Ne}
\newcommand {\thstar} {\theta_\ast}
\newcommand {\Lie} [1] {{\cal L}_{#1}}
\newcommand {\gbar} {\overline g}
\newcommand {\ddeps} 
	{\left. {\frac d {d\epsilon}} \right|_{\epsilon = 0}}
\newcommand {\DDomega} [2] {\D {#1} \D {#2} \omega}
\newcommand {\Abar} {\overline A}
\newcommand {\Omegabar} {\overline \Omega}
\newcommand {\Omegabbar} {\overline {\overline \Omega}}
\newcommand {\Phibar} {\overline \Phi}
\newcommand {\Psibar} {\overline \Psi}
\newcommand {\Abarbar} {\overline {\overline A}}
\newcommand {\Omegabarbar} {\overline {\overline \Omega}}
\newcommand {\tomega} {\tilde {\omega}}
\newcommand {\tW} {\tilde {W}}

\newcommand {\BM} {basilar membrane }
\newcommand {\micron} {$\mu$m}
\begin{document}
\title{Modeling Elastic Shells Immersed in Fluid}
\thanks{
Computation was performed at the Pittsburgh
Supercomputing Center under
allocation MCA93S004 from the MetaCenter Allocations Committee.
}
% \thanks{...}% At most 5 thanks
%
\author{Edward Givelberg}
%% \address{Department of Mathematics,
%% University of Michigan, 
%% 525 East University Avenue, 
%% Ann Arbor, MI 48109
%% }
% \author{...}\address{...}
% \author{...}\address{...}
%
% \date{September 1, 2001}
%
\begin{abstract} 
%================================================================
% \input {abstract}
We describe a numerical method to simulate 
an elastic shell immersed in a viscous incompressible fluid.
The method is developed as an extension of the immersed boundary method
using shell equations based on the
Kirchhoff-Love and the planar stress hypotheses.
A detailed derivation 
of the shell equations used in the numerical
method is presented.
This derivation 
as well as the
numerical method, use techniques of differential geometry
in an essential way.
Our main motivation for the development of this method
is its use in the construction of 
a comprehensive 
three-dimensional computational model of the cochlea
(the inner ear).
The central object of study within the cochlea is 
the ``basilar membrane'', 
which is immersed in fluid
and
whose elastic properties rather resemble those of a shell.
We apply the method to a specific example, which is a prototype of
a piece of the basilar membrane and study the convergence of the
method in this case.
Some typical features of cochlear mechanics are already captured
in this simple model.
In particular, 
numerical experiments have shown a traveling wave propagating from
the base to the apex of the model shell in response to external
excitation in the fluid.
\end{abstract}
%================================================================
%
%
%% \subjclass{74F10, 74K25, 74B05, 76Z99, 92C05}
%
%% \keywords{
%% elastic shell theory, 
%% Navier-Stokes equations,
%% immersed boundary method,
%% cochlea
%% }
%
\maketitle
%%
%%
%%
%================================================================
% \input {intro}
\section {Introduction}
\label {intro}

\newcommand {\Bekesy} {B\'{e}k\'{e}sy}

This paper describes a general method of simulation of
an elastic shell 
immersed in a viscous incompressible fluid.
The method is developed as an extension of the immersed boundary method
originally introduced by Peskin and McQueen to study the blood flow 
in the heart.
The immersed boundary method has proved to be particularly useful for
computer simulation of various biofluid dynamic systems.
In this framework the elastic (and possibly active) biological tissue
is treated as a collection of elastic fibers immersed in 
a viscous incompressible fluid.
This formulation of the method together with references 
to many applications
can be found in \cite{Peskin1994}.
A partial list of the applications of the immersed boundary method
includes 
in addition to the blood flow in the heart
(see the extensive work of Peskin and McQueen, e. g.
\cite{Peskin1997, McQueen2001})
also
platelet aggregation during blood clotting
\cite{Fauci1993},
flow of suspensions \cite{Fogelson1988, Sulsky1991},
aquatic animal locomotion \cite{Fauci1993},
a two-dimensional model of cochlear fluid mechanics
\cite{Beyer1992}
and
flow in collapsible tubes
% \cite{Rosar1994}.
\cite{Rosar2001}.
For a recent review of the immersed boundary method, see
\cite{Peskin2002}.

Many man-made materials can be modeled as elastic shells
and elastic shells are also ubiquitous in nature,
however
our motivation for the
development of the method comes from the study of the cochlea.
The auditory signal processing in the cochlea depends crucially
on the dynamics of the basilar membrane,
which is immersed in a viscous 
incompressible fluid of the cochlea,
and despite its name the basilar membrane is actually an elastic shell.
The numerical method for elastic shell--fluid interaction presented
here
was subsequently used in the
construction of a complete three-dimensional computational model
of the macro-mechanics of the cochlea which incorporates the intricate
curved cochlear anatomy. 
The results of this work will be reported in future publications.

The cochlea is the part of the inner ear where sound waves
are transformed into electrical pulses which are 
carried by neurons to the brain.
It is a small snail-shell-like cavity in the temporal bone,
which has two openings, the oval window and the round window.
The cavity is filled with fluid, which is sealed in 
by two elastic membranes covering 
these windows.
The spiral canal of the cochlea
is divided lengthwise by the long and narrow basilar membrane 
into two passages
that connect with each other at the apex.
External sounds set the ear drum in motion,
which is conveyed to the inner ear by
three small bones of the middle ear.
These bones function as an impedance matching device,
focusing the energy of the ear drum on the oval window of the
cochlea.
This piston-like motion against the oval window
displaces the fluid of the cochlea generating traveling waves
that propagate along the basilar membrane.
The vibrations of the basilar membrane are detected by
thousands of microscopic sensory receptors,
called hair cells,
located on the surface of the basilar membrane.
The auditory signal processing in the cochlea is completed
by the hair cells converting these mechanical stimuli 
into action potentials in the neurons attached to them,
relaying this information to the brain.

Practically everything we know about
the passive wave propagation in the cochlea was discovered in the 1940s 
by Georg von \Bekesy\,
% \cite{vonBekesy},
who carried out experiments in cochleae extracted 
from human cadavers.
Von \Bekesy\ observed that 
a pure tone input sound
generates a traveling wave which reaches its peak at 
a speciffic location along the \BM\,
exciting only a narrow band of hair cells.
This characteristic location depends on the tone's frequency.
By pinching the basilar membrane with a tiny probe and observing the
resulting displacement von \Bekesy\ discovered that 
the basilar membrane is,
in fact, not a membrane, i. e. it is not under inner tension,
but an elastic shell, whose compliance varies
exponentially along the membrane.
Von \Bekesy\'s extensive experimental work is summarized in
his book "Experiments in Hearing"
\cite{vonBekesy}.
% varying as
% $e^{-\lambda x}$ along the basilar membrane, with
% $\lambda^{-1} = 0.7 \text{cm}$.
% The motion of the basilar membrane
% is therefore crucial to
% understanding the auditory signal processing in the cochlea.
For an excellent summary of more recent work on the cochlea, see
\cite{Geisler1998}.

Cochlear mechanics has been an active area of research ever since
von \Bekesy's fundamental contributions, yet
many important questions are still open.
Presently there is no complete understanding of 
the mechanisms responsible for the extreme sensitivity, 
sharp frequency
selectivity and broad dynamic range of hearing.
% Mathematicians have also contributed to the understanding of the
% cochlear mechanics.
% The 197? paper of Neu and Keller [] ....
The most rigorous mathematical analysis of the cochlea
was carried out by
Leveque, Peskin and Lax in \cite{Leveque1988}.
In their model the cochlea is represented by
a two-dimensional plane 
(i.e. a strip of infinite length and infinite depth)
and the basilar membrane, by a straight line of harmonic oscillators
dividing the fluid plane into two halves.
The linearized equations are reduced to a functional equation by
applying the Fourier transform in the direction parallel to the
basilar membrane and then solving the resulting ordinary differential
equations in the normal direction.  The functional equation derived in
this way is solved analytically, and the solution is evaluated both 
numerically and also asymptotically (by the method of stationary phase).
This analysis reveals that the waves in the cochlea resemble
shallow water waves, i.e. ripples on the surface of a pond.
A distinctive feature
of this paper is the (then speculative) consideration 
of negative basilar 
membrane friction, i.e., of an amplification mechanism operating within
the cochlea.

Like the cochlea, even the simplest
three-dimensional fluid--shell systems appear to be too
difficult to analyze,
but using the methods presented here it is possible to construct
computational models for such systems.
The construction of computational models for fluid--shell
configurations is however not easy,
and the development of these methods
required supercomputing resources.
Nevertheless, with the rapid advances in computer technology,
such computations will soon be feasible on a workstation.
The cochlea example is important also because
the fluid--shell system simulation in this case can be compared
with
the extensive body of theoretical and experimental research.

The theory of elastic plates and shells is a classical mathematical
subject, which is also an active area of contemporary research
(see for example
\cite{Ciarlet, Niordson, Vorovich}).
For completeness,
we derive 
in section \ref {ChapterBundle} 
the elastic shell equations upon which
the numerical method is based.
Shell theory is very naturally described in the language of
differential geometry and it is perhaps in this respect that
the presentation in section \ref {ChapterBundle}
somewhat differs from other accounts of the subject.
Imagine a material composed of very rigid straight line segments 
(fibers) coupled together and assume that these fibers are 
perpendicular to some imaginary surface $S$ 
in the middle of the material.
In other words, the material has a structure of a normal vector bundle
over the surface $S$.
We shall assume that the base surface $S$ is free to undergo 
arbitrary small elastic deformations. 
The Kirchhoff-Love hypothesis, 
% which occupies a central place in shell theory, 
is the assumption that the deformation
of the bundle is such that the fibers of the deformed bundle are
perpendicular to the deformed middle surface,
and that these fibers
are not stretched during the deformation. 
In the language of differential geometry 
this means
that the shell has the structure of a normal vector bundle 
which is preserved under elastic deformations. 
We shall call such a material an 
{\em elastic bundle}.
The Kirchhoff-Love hypothesis implies that a deformation
of the base surface completely determines the deformation of 
the whole bundle.
To ensure that the fibers do not intersect each other we must
assume that their length, i.e., the thickness of the bundle,
is smaller than the radius of maximal curvature
of the base surface.

Taking the three-dimensional linear theory of elasticity as our
starting point, and using the Kirchhoff-Love hypothesis, 
% we derive
it is possible to derive
the equations which completely determine small deformations of
elastic bundles. 
It turns out however that a realistic model has to satisfy
the hypothesis of planar stress.
Strictly speaking, in linear elasticity the planar stress hypothesis
is not consistent with the Kirchhoff-Love assumption.
We utilize a common approach to ``reconciling'' the two assumptions
(e.g see \cite{Hughes2000}).
% These equations constitute a system of three partial differential
% equations in three unknown functions.
% We would like to emphasize that 
% The equations are exact and do not
% contain thin shell approximations. 
%% The only assumption we make on the
%% thickness of the bundle is the one related to the curvature radius
%% stated above.
The resulting equations
constitute a system of three partial differential
equations in three unknown functions.
They
express the elastic force exerted by the shell 
% displaced from equilibriuim 
as a linear fourth order differential
operator applied to the displacement vector field.
This differential operator is
intrinsic to the base surface of the bundle
and its coefficients are tensorial quantities determined by
the geometry and the elastic properties of the bundle.
% We can subsequently simplify these
% equations by
% % derive shell theory equations from 
% % the bundle equations by
% assuming that the thickness of the bundle is small compared to its
% other dimensions.
This intrinsic geometric formulation is 
% at the basis of the
essential in the formulation of the
numerical method described in the following section.

Section \ref {ChapterIBM} outlines 
the immersed boundary method for shells. 
This is a modification of the immersed boundary
method as described in \cite{Peskin1994}. 
The main difference in the algorithm is
in the computation of the force that the material applies to the fluid. 
Here it is computed by discretizing the shell equations 
described in section \ref {ChapterBundle}.
We describe in detail the discrete differential operators
defined on the shell surface which are used in the force computation.

%
% Accordingly, the example we present in section \ref {ChapterApplication}
In section \ref {ChapterApplication}
we describe a test model in which the shell
resembles a piece of a basilar membrane immersed in fluid.
We study the convergence of the algorithm in this case
and describe the numerical experiments carried out with this model.
These experiments have reproduced some of the typical features of 
cochlear mechanics, such as the traveling wave propagating along
the basilar membrane in response to external excitation of the fluid.

% The main difficulty in constructing this model has been the
% large scale of the immersed boundary computations required.
% In the last section we discuss the complexity of immersed
% boundary computations
% and
% some of the challenges in
% the construction of a computational model 
% of the cochlea.
%================================================================

% \input {acknowledge}

This work is based on the author's Ph.D. thesis completed at the 
Courant Institute (NYU) under the supervision of 
Professor Charles S. Peskin. I would like to thank Professor Peskin
for his support, encouragement and patient guidance.
I would also like to thank Professor David McQueen 
for many conversations
in which I learnt a lot about scientific computing.
I'd like to thank Professor Karl Grosh for helping me understand shell
theory.
%================================================================
%
%================================================================
% \input {bundle}
`

\section {Linear Elastic Deformations of Vector Bundles (Shell Theory)}
\label {ChapterBundle}

Let $ M \subset \Rthree $
be a normal bundle over a surface $S$.
% be a bundle of line segments in the Euclidean space
% with the surface $S$ its base. 
We describe $S$ by a coordinate chart 
$ \phi : \Omega \to \Rthree $,
where $\Omega$ is a domain in  $\Rtwo$.
Let 
$ n : \Omega \to M $ be the unit normal vector field on $S$.
The natural chart for $M$ is
$ \Phi : \Omega \times (-\thickness, \thickness) \to M $
given by
$$
\Phi(q_1, q_2, t) = \phi(q_1, q_2) + t n(q_1, q_2),
$$
and we assume, for simplicity, that the thickness of the bundle
$(= 2\thickness)$ is constant.
Let 
$ \Phi_\epsilon : \Omega \times (-\thickness, \thickness) \to \Rthree $
describe a 1-parameter family of deformations of $M$,
such that
$ \Phi_0 = \Phi $.
Our basic assumption is that each $\Phi_\epsilon$ preserves the bundle
structure of $M$.
Clearly, such a deformation is completely determined by the
deformation of the base space $S$. 
In the framework of linear elasticity we work with infinitesimal
deformations, i.e., vector fields. 
We will assume that $S$ is free to undergo an arbitrary
infinitesimal deformation
$ \psi = \frac d {d\epsilon}|_{\epsilon = 0} \pe $
and we shall determine the corresponding infinitesimal deformation
$ V = \frac d {d\epsilon}|_{\epsilon = 0} \Pe $
of $M$.

We begin with some preliminary remarks about the geometry of $M$.
Throughout this chapter the indices are raised and lowered with
respect to $h$, the metric of $S$.
We adopt the convention that greek indices take the values $1, 2$,
while roman indices run
through $1, 2, 3$.
The standard metric on $\Rthree$ will be denoted by $\delta$:
$$
\delta(X, Y) = <X, Y> = X \cdot Y,
$$
$ \partial_\alpha = \frac \partial {\partial q_\alpha} $ is a partial
derivative and
$\overline\nabla$ is the standard flat connection of $\Rthree$.

%---------------------------------------------------------------
\subsection {The Geometry of the Vector Bundle}

Let
$ h = \phi^\ast \left( \delta|_S \right) $
denote the metric of $S$ in the chart $\phi$. 
The components of $h$ are
$$
h_{\alpha\beta} = \dphi {\alpha} \cdot \dphi {\beta}.
$$
Throughout this paper greek tensor indices are raised and lowered
with respect to the metric $h$.
Let $N_p$ denote the unit normal vector to $S$ at the point
$ p \in M $.
Thus, $ n(q_1, q_2) = N_{\phi(q_1, q_2)} $.
We shall often abuse notation and identify $n$ with $N$.
Similarly for other vector fields.
The second fundamental form of $S$ is the symmetric bilinear form
acting on vectors tangent to the surface $S$ defined by
$$
b(X, Y) = <\DD X N, Y>.
$$
Its components in the chart $\phi$ are
$$
b_{\alpha\beta} = \dn {\alpha} \cdot \dphi {\beta}.
$$
The central role in the following analysis is played by the map
$\theta = \theta_t$, which we define to be the flow of $N$:
$$
\theta_t(p) = p + t \, N_p, \ \ \ \ \, p \in S.
$$
The differential of $\theta$ is the map
$$
\thstar : T_pS \to T_{\theta(p)}S_t,
$$
where $ S_t = \theta_t(S) $.
Let $ X \in T_pS $ and let $\sigma$ be a curve on $S$ such that
$ \sigma(0) = p $, 
$ \sigma'(0) = X $.
Then
\begin {eqnarray*}
	\thstar X & = & \thstar \sigma' \\
			& = & (\theta \circ \sigma)' \\
			& = & (\sigma + t \, N \circ \sigma)' \\
			& = & X + t \, \DD X N
\end {eqnarray*}
We shall regard $\thstar$ as a symmetric tensor field on $S$.
The components $\theta_{\alpha\beta}$
of $\thstar$ can be expressed in the chart $\phi$
as follows:
\begin {eqnarray}
\theta_{\alpha\beta} & = & < \thstar(\dphi \alpha), \, \dphi \beta> 
						\nonumber \\
	& = & < \dphi \alpha + t \, \DD {\dphi \alpha} N, \, \dphi \beta>
						\nonumber \\
	& = & \dphi \alpha \cdot \dphi \beta 
		+ t \, \dn \alpha \cdot \dphi \beta
						\nonumber \\
	& = & h_{\alpha\beta} + t \, b_{\alpha\beta}.
\end {eqnarray}
An important observation which will be useful later is:
\begin{equation}
\DD X N = \DD {\thstar X} N
\end{equation}
for any tangent vector $X$.
Indeed, let 
$ \tilde{\sigma} = \theta \circ \sigma $.
Then
$ \thstar X = \tilde{\sigma}' $
and
\begin {eqnarray*}
	\DD {\thstar X} N & = & \DD {\tilde {\sigma}'} N \\
		& = & \lim_{\epsilon \to 0} \frac 
			{N(\tilde{\sigma}(\epsilon)) - N(\tilde{\sigma}(0))}
			\epsilon \\
		& = & \lim_{\epsilon \to 0} \frac 
			{N(\sigma(\epsilon)) - N(\sigma(0))}
			\epsilon \\
		& = & \DD X N
\end {eqnarray*}
%
%
%
% The metric of $M$ is the induced metric which we denote by
% $ \overline g $.
We denote the metric of $M$ in the chart $\Phi$ by $\ggg$.
Its components,
$
\ggg_{ij} = \dPhi i \cdot \dPhi j
$,
are
%
% We have therefore
% \begin {eqnarray*}
% \dPhi 3 & = & n \\
% \dPhi \alpha & = & \theta_\alpha{}^\sigma \dphi \sigma
% \end {eqnarray*}
% and
\begin {eqnarray*}
\ggg_{\alpha\beta} & = & \theta_\alpha{}^\sigma \theta_{\sigma\beta}
\\
\ggg_{3\alpha} & = & 0
\\
\ggg_{33} & = & 1
\end {eqnarray*}
The collection of parallel surfaces $S_t$ forms a foliation of
$M$. 
This foliation carries the induced connection $\nabla$ given by
$$
\D X Y = \DD X Y - <\DD X Y, \, N>N
$$
for any two vector fields $X$ and $Y$ which are tangent to the
foliation.
The second fundamental form of this foliation is the symmetric
bilinear form acting on tangent vectors defined by
$$
B(X, Y) = <\DD X N, \, Y>.
$$
% where X and Y are vector fields tangent to the foliation.
Thus
$$
\D X Y = \DD X Y + B(X, Y) N.
$$
The components of $B$ in the chart $\Phi$ are
\begin {eqnarray}
\label{eq: Bab}
B_{\alpha\beta} & = & 
		(\Phi^\ast B)(\partial_\alpha, \partial_\beta) 
										\nonumber \\
	& = & B(\dPhi \alpha, \dPhi \beta) 
										\nonumber \\
	& = & <\DD {\dPhi \alpha} N, \dPhi \beta> 
										\nonumber \\
	& = & <\DD {\thstar(\dphi \alpha)} N, \dPhi \beta> 
										\nonumber \\
	& = & <\DD {\dphi \alpha} N, \dPhi \beta> 
										\nonumber \\
	& = & \dn \alpha \cdot (\theta_\beta{}^\sigma \dphi \sigma) 
										\nonumber \\
	& = & b_{\alpha\sigma} \theta_\beta{}^\sigma  
										\nonumber \\
	& = & \theta_\alpha{}^\sigma b_{\sigma\beta}.
\end {eqnarray}
The map $\theta$ plays an important role because it allows us
to extend any tensor $A$ on the base surface $S$ to a tensor
$\thstar A$ on the whole bundle $M$.

We conclude the description of the geometry of $M$ with its volume
element $dv$. 
It can be decomposed as follows
$$
dv = d A_t \, dt,
$$
where $dA_t$ is the area element of $S_t$.
Since $ S_t = \theta(S) $, we have
$$
dA_t = \det(\thstar) \, dA_0.
$$
Notice than in $\Rthree$ 
\begin {eqnarray*}
\det(\thstar) & = & \det(h + tb) \\
	& = & 1 + (\tr \,b) \, t + (\det \, b) \, t^2 \\
	& = & 1 + H \, t + K \, t^2,
\end {eqnarray*}
where $H$ and $K$ denote the mean curvature and the Gaussian curvature
of the surface $S$ respectively.

%---------------------------------------------------------------
\subsection {The Infinitesimal Deformation of an Elastic Bundle}

We shall assume that the deformation field of the base
surface $S$ is given by
$$
\psi = \omega \, N + W ,
$$
where $ W = W^\sigma \dphi \sigma $ 
is an arbitrary vector field tangent to $S$
and $\omega$ is an arbitrary function on $S$.
We extend $\omega$ to the whole bundle $M$ by the normal
flow, but we'll continue to write $\omega$ instead of 
$\thstar \omega = \omega \circ \theta$.
In this section we will show that the 
corresponding deformation vector field
of the elastic bundle $M$ is given by
\begin {equation}
\label {eq: DeformationVF}
V = \omega \, N + \thstar W - t \, \Tom
\, ,
\end {equation}
where 
$$
 \Tom(q_1, q_2, t) = \Tom(q_1, q_2) 
	= h^{\mu\nu} \, \domega \mu \, \dphi \nu .
$$
From the assumption that the deformations should preserve the
bundle structure it follows that $\Pe$ has the following form:
$$
\Pe = \pe + t \, \Ne
\, ,
$$
where $\pe$ is a coordinate system for the deformed base space
and $\Ne$ is the normal to it.
We have 
$ \phi_0 = \phi $,
$ n_0 = n $
and
\begin {eqnarray}
\label {eq: V}
V & = & \ddeps \Pe 
			\nonumber \\
	& = & \psi + t \, \ddeps \Ne
\, .
\end {eqnarray}
On one hand
\begin{equation}
\label {eq: dNn}
\left( \ddeps \Ne \right) \cdot n 
	= \frac 1 2 \, \ddeps (\Ne \cdot \Ne) = 0
\, ,
\end{equation}
while on the other hand
\begin {eqnarray}
\label {eq: dNtang}
\left( \ddeps \Ne \right) \cdot \dphi \alpha 
	& = & \ddeps (\Ne \cdot \dpe \alpha)
		- n \cdot \ddeps \dpe \alpha 
												\nonumber \\
	& = & - n \cdot \dpsi \alpha 
												\nonumber \\
	& = & - \partial_\alpha (\psi \cdot n) + \psi \cdot \dn \alpha 
												\nonumber \\
	& = & - \domega \alpha + \psi \cdot b_\alpha{}^\sigma \dphi \sigma 
												\nonumber \\
	& = & - \domega \alpha + b_\alpha{}^\sigma W_\sigma
\, .
\end {eqnarray}
Using (\ref{eq: dNn}) and (\ref{eq: dNtang}) in (\ref{eq: V})
we have
\begin {eqnarray*}
V & = & \omega \, N + W + t \, 
		(- \domega \alpha  + b_\alpha{}^\sigma \, W_\sigma)
		\, h^{\alpha\tau} \, \dphi \tau \\
	& = & \omega \, N + (h_\alpha{}^\sigma + t \, b_\alpha{}^\sigma) 
		\, W_\sigma \, h^{\alpha\tau} \, \dphi \tau
		- t \, \domega \alpha \, h^{\alpha\tau} \, \dphi \tau \\
	& = & \omega \, N + \theta_\alpha{}^\sigma \, W_\sigma 
		\, h^{\alpha\tau} \, \dphi \tau
		- t \, \domega \alpha \, h^{\alpha\tau} \, \dphi \tau \\
	& = & \omega \, N + \thstar W - t \, \Tom
\, ,
\end {eqnarray*}
which is what was claimed in (\ref {eq: DeformationVF}).
%
% and
% $$
% \Ne = \frac {\dpe {1} \times \dpe {2}} {|\dpe {1} \times \dpe {2}|}
% $$
% 
% 

%---------------------------------------------------------------
\subsection {The Strain Tensor}

Let $ \gbar = \delta |_M $ be the metric of the bundle $M$.
The strain tensor corresponding to the deformation vector field
$V$ is the symmetric tensor $e$ defined by
$$
e = \frac 1 2 \Lie V \gbar,
$$
the Lie derivative of $\gbar$ in the direction of $V$.
We now show that for any tangent vector field $X$, 
\begin{equation}
\label{eq:ZeroStrain}
e(X, N) = 0.
\end{equation}
Let 
$ g_\epsilon = \Phi^\ast_\epsilon (\gbar) $
be the metric of the deformed bundle in the chart $\Pe$.
Its components are:
$$
(g_\epsilon)_{ij} = \dPe i \cdot \dPe j.
$$
The components of the strain tensor in the chart $\Phi$ are
$$
e_{ij} = \frac 1 2 \ddeps
		\left( \dPe i \cdot \dPe j \right)
$$
We have
\begin {eqnarray}
\label {eq: StrainTensor33}
e_{33} & = & \dPhi 3 \cdot \ddeps \dPe 3 
											\nonumber \\
	& = & n \cdot \ddeps n_\epsilon 
											\nonumber \\
	& = & 0
\end {eqnarray}
and
\begin {eqnarray}
\label {eq: StrainTensor3T}
e_{3\alpha} & = & \frac 1 2 \ddeps
	\left( n_\epsilon \cdot \dpe \alpha 
			+ t \, \Ne \cdot \dne \alpha \right) 
											\nonumber \\
		& = & \frac 1 2 \left( - (\dpsi \tau \cdot n) h^{\tau\sigma}
			\dphi \sigma \cdot \dphi \alpha
			+ n \cdot \dpsi \alpha \right) 
											\nonumber \\
		& = & 0,
\end {eqnarray}
which shows (\ref{eq:ZeroStrain}).
We can now express the strain tensor in terms of the function
$\omega$ and the vector field $W$:
\begin {eqnarray*}
e & = & \frac 1 2 \Lie V \, \gbar \\
	& = & \frac 1 2 \Lie {\omega N + \thstar W - t \, \Tom} \, \gbar \\
	& = & \frac 1 2 \Lie {\omega N} \, \gbar 
			+ \widehat {\nabla (\thstar W)}
			- t \, \widehat {\nabla \Tom}.
\end {eqnarray*}
Here $ \hat {A} $  stands for the symmetrization of A:
$$
\hat {A}_{\mu\nu} = \frac 1 2 (A_{\mu\nu} + A_{\nu\mu}).
$$
For tangential fields $X$, $Y$ we have
$$
\Lie {\omega N} \, \gbar(X, Y) = \omega \, N <X, Y>
	- <\Lie {\omega N} X, Y>
	- <X, \Lie {\omega N} Y> .
$$
Since
\begin {eqnarray*}
	\Lie {\omega N} X & = & \omega N X - X(\omega N) \\
		& = & \omega \Lie N X - (X\omega) N
\end {eqnarray*}
we have
$$
<\Lie {\omega N} X, Y> = \omega <\Lie N X, Y>
$$
and
$$
\Lie {\omega N} \gbar (X, Y) = \omega \Lie N \gbar(X, Y).
$$
Finally the strain tensor is given by
\begin {equation}
\label {eq: StrainTensor}
e = \omega B + \widehat {\nabla (\thstar W)}
			- t \, \widehat {\nabla \Tom}.
\end {equation}

For future reference we obtain the expressions for the
components of the tensors
$\nabla \Tom$ and $\nabla (\thstar W)$ in the chart $\Phi$:
\begin {eqnarray}
\label {eq: DTom}
\nabla_\alpha \Tom_\beta & = & \nabla \Tom(\dPhi \alpha, \, \dPhi \beta)
                                                    \nonumber \\ 
    & = & \theta_\alpha{}^\sigma \theta_\beta{}^\tau
        \nabla \Tom(\dphi \sigma, \, \dphi \tau)
                                                    \nonumber \\
    & = & \theta_\alpha{}^\sigma \theta_\beta{}^\tau
        <\D {\dphi \sigma} \Tom, \, \dphi \tau>
                                                    \nonumber \\
    & = & \theta_\alpha{}^\sigma \theta_\beta{}^\tau
        <\DD {\dphi \sigma} \Tom, \, \dphi \tau>
                                                    \nonumber \\
    & = & \theta_\alpha{}^\sigma \theta_\beta{}^\tau
        \partial_\sigma \Tom \cdot \dphi \tau
                                                    \nonumber \\
    & = & \theta_\alpha{}^\sigma \theta_\beta{}^\tau
        \DDomega \sigma \tau,
\end {eqnarray} 
where
$ \DDomega \sigma \tau $
are the components of the Hessian of $\omega$ on $S$
in the chart $\phi$,
and
\begin {eqnarray}
\label {eq: DthW}
\nabla(\thstar W)_{\alpha\beta} & = &
        <\D {\dPhi \alpha} (\thstar W), \, \dPhi \beta>
                                                    \nonumber \\
    & = & \theta_\alpha{}^\sigma \theta_\beta{}^\tau
        <\D {\dphi \sigma} (\thstar W), \, \dphi \tau>
                                                    \nonumber \\
    & = & \theta_\alpha{}^\sigma \theta_\beta{}^\tau
        <\D {\dphi \sigma} (W^\mu \theta_\mu{}^\nu \dphi \nu),
            \, \dphi \tau>
                                                    \nonumber \\
    & = & \theta_\alpha{}^\sigma \theta_\beta{}^\tau
        \D \sigma (W^\mu \theta_{\mu\tau})
                                                    \nonumber \\
    & = & \theta_\alpha{}^\sigma \theta_\beta{}^\tau
%           \theta_\tau{}^\mu (\D \sigma W_\mu)
            \theta_{\tau\mu} (\D \sigma W^\mu)
        + \theta_\alpha{}^\sigma \theta_\beta{}^\tau
            (\D \sigma \theta_{\mu\tau}) W^\mu
\end {eqnarray}
Substituting
(\ref {eq: DTom}) 
and 
(\ref{eq: DthW}) 
in (\ref {eq: StrainTensor})
we obtain
\begin{equation}
\label{eqn:straincomponents}
e_{\alpha\beta} =
B_{\alpha\beta} 
\,
\omega
\, - \,
t \, \theta_\alpha{}^\sigma \theta_\beta{}^\tau
\, \DDomega \sigma \tau 
\, +  \,
\theta_\alpha{}^\sigma \ggg_\beta{}^\tau
\, (\D \sigma W_\tau)
\, + \,
\theta_\alpha{}^\sigma \theta_\beta{}^\tau
(\D \sigma \theta_\tau{}^\mu) 
\, W_\mu
\, .
\end{equation}

%================================================================
% \input {planestress}

\subsection {The Plane Stress Assumption.}

In linear elasticity the stress tensor $\sigma$ is related to the strain
tensor $e$ by the generalized Hooke's law:
$$
\sigma^{ij} = C^{ijkl} e_{kl}
$$
where 
$$
C_{ijkl} = C_{ijkl}(q_1, q_2, t)
$$
denote the components of the elasticity 4-tensor
$C$ in the chart $\Phi$.
For a homogeneous and isotropic material
% We shall assume 
the elasticity tensor is of the form
\begin{equation}
C^{ijkl} = \lambda \, \ggg^{ij} \ggg^{kl} 
			+ \mu \, \ggg^{ik} \ggg^{jl} + \mu \, \ggg^{il} \ggg^{jk}
\, .
\end{equation}
% (see Marsden....).
In this case Hooke's law of linear elasticity takes the form
\begin{equation}
\label {eqn:hookeslaw}
\sigma_{ij} = \lambda \, \ggg^{mn} e_{mn} \, \ggg_{ij} 
				+ 2 \mu \, e_{ij}
			= \lambda \, \tr(e) \, \ggg_{ij} + 2 \mu \, e_{ij} \, .
\end{equation}
The strain energy corresponding to the strain $e$ is
defined by
$$
{\cal E} = \frac 1 2 \int_M C(e, e) \, dv
=  \frac 1 2 \int_M C^{ijkl} e_{ij} e_{kl}
$$
%
%  The tensor $C$ has the following symmetries
%  $$
%  C_{ijkl} = C_{klij}
%  $$
%  $$
%  C_{ijkl} = C_{jikl}
%  $$
%  and we think of it as a symmetric bilinear form acting on 2-tensors
%  with the additional property that
%  $$
%  C(A, B) = C(\hat{A}, B)
%  $$
%  for any 2-tensors $A$, $B$.
%
%
% We now develop the shell equations using the plane stress assumption.
%
It is now possible to derive a shell theory based solely on the
Kirchhoff-Love hypothesis using the expression for the strain tensor
(\ref{eqn:straincomponents}).
It turns out however that a more realistic shell model is obtained
using 
the assumption of plane stress: 
\begin{equation}
\label {eqn:planestress}
\sigma_{33} = 0.
\end{equation}
In linear elasticity this assumption contradicts the
Kirchhoff-Love hypothesis.
We follow a standard approach in modifying the 
Kirchhoff-Love hypothesis
(see \cite{Hughes2000}):
% We therefore modify the assumption on the strain tensor as follows:
we no longer assume that
$$
e_{33} = 0,
$$
but we continue to assume that
$$
e_{3\alpha} = 0.
$$
\newcommand {\ebar} {\overline e}
Let
$\ebar$ denote the tangential part of the strain $e$.
From (\ref{eqn:hookeslaw}) and (\ref{eqn:planestress})
it follows that 
$$
0 = \sigma_{33} = \lambda \ggg^{mn} e_{mn} + 2 \mu e_{33} 
= \lambda \ggg^{\mu\nu} e_{\mu\nu} + (\lambda + 2 \mu) e_{33}
$$
and therefore
$$
e_{33} = - \frac \lambda {\lambda + 2 \mu} \ggg^{\mu\nu} e_{\mu\nu}
= - \frac \lambda {\lambda + 2 \mu} tr \ebar
\, .
$$
We use the last equation
to simplify the expression for the
strain energy:
$$
C(e, e) = \lambda \, (g^{mn} e_{mn})^2 
			+ 2 \mu \, g^{mn} g^{kj} \, e_{mj} e_{kn}
		= \lambda (\tr e)^2 + 2\mu \tr (e^2)
\, .
$$
Since
$$
\tr e = \tr \ebar + e_{33} 
= (1 - \frac \lambda {\lambda + 2 \mu})~\tr \ebar
= \frac {2 \mu} {\lambda + 2 \mu}~\tr \ebar
$$
%
%%  Since for any $2 \times 2$ symmetric tensor $A$
%%  $$
%%  \tr(A^2) = (\tr A)^2 + 2 \det A
%%  $$
and
$$
\tr(e^2) = 
	\tr(\ebar^2) + e_{33}^2
	= \tr(\ebar^2)
+ \frac {\lambda^2}  {(\lambda + 2 \mu)^2} (\tr \ebar)^2
$$
we obtain
\begin{eqnarray*}
C(e, e) & = & 
	4 \lambda \frac {\mu^2} {(\lambda + 2 \mu)^2} (\tr \ebar)^2
	+ 2 \mu \frac {\lambda^2}  {(\lambda + 2 \mu)^2} (\tr \ebar)^2
	+ 2 \mu \tr(\ebar^2)
%% + 2 \mu (\tr \ebar)^2
%% + 4 \mu \det(\ebar)
\\
& = & 
	\frac {2 \lambda \mu} {\lambda + 2 \mu} (\tr \ebar)^2
	+ 2 \mu \tr (\ebar^2)
\, .
\end{eqnarray*}

\subsection {Variation of the Strain Energy} 
\label {sect:VarStrainE}

In order to obtain the equilibrium equations we proceed to
calculate the variation of the strain energy. 
Let $\tomega$ be an arbitrary function and $\tW$ an arbitrary
vector field on $S$. 
\newcommand {\tV} {\tilde V}
\newcommand {\te} {\tilde e}
The corresponding variation of the energy is
$$
\delta {\cal E} = \lim_{\epsilon \to 0}
    \frac {{\cal E}(\omega + \epsilon \tilde{\omega},
            W + \epsilon \tilde{W}) - {\cal E}(\omega, W)} \epsilon
\, .
$$
Let
$\te = \ebar(\tV)$
be the strain corresponding to the deformation $\tV$.
Since
% the displacement field $V$ determines the strain
% tensor $\ebar$
$$
\ebar(V + \tV) = \ebar(V) + \ebar(\tV)
$$
the variation of the strain energy is
\begin{eqnarray*}
\delta {\cal E}
& = &
\frac {\lambda \mu} {\lambda + 2 \mu} \,
\int_M 
	\tr \ebar~ \tr \te
\, dv
+ \mu \int_M 
	\tr (\ebar \, \te)
\, dv
\\
& = &
\int_\Omega
\int_{-h_0}^{h_0}
	\left(
		\frac {\lambda \mu} {\lambda + 2 \mu} \,
		\ggg^{\alpha\beta} \ggg^{\gamma\delta}
		+ \mu \,
		\ggg^{\alpha\gamma} \ggg^{\beta\delta}
	\right)
	e_{\alpha\beta} \te_{\gamma\delta}
	\, \det (\thstar) 
\, dt 
\, dA_0
\\
& = &
\int_\Omega
\int_{-h_0}^{h_0}
	\Lambda^{\alpha\beta\gamma\delta} 
	e_{\alpha\beta} \te_{\gamma\delta}
\, dt 
\, dA_0,
\end{eqnarray*}
where we define
$$
\Lambda^{\alpha\beta\gamma\delta}(q_1, q_2, t) =
\left(
\frac {\lambda \mu} {\lambda + 2 \mu} \,
\ggg^{\alpha\beta} \ggg^{\gamma\delta}
+ \mu \,
\ggg^{\alpha\gamma} \ggg^{\beta\delta}
\right)
\, \det (\thstar)
\, .
$$
For the calculations that follow it is useful to note the following
symmetries  of $\Lambda$:
$$
\Lambda^{\alpha\beta\gamma\delta}
=
\Lambda^{\beta\alpha\gamma\delta}
=
\Lambda^{\alpha\beta\delta\gamma}
=
\Lambda^{\gamma\delta\alpha\beta}
\, .
$$
We use (\ref {eqn:straincomponents})
and integrate by parts omitting the boundary terms
to rewrite
the energy variation as follows
\begin{eqnarray*}
\delta {\cal E}
& = &
\int_\Omega
\int_{-h_0}^{h_0}
\Lambda^{\alpha\beta\gamma\delta} e_{\alpha\beta} 
\left(
	B_{\gamma\delta} \tomega
	- t \, \theta_\gamma{}^\sigma \theta_\delta{}^\tau
	\D \sigma \D \tau \tomega
	+ \theta_\gamma{}^\sigma \ggg_\delta{}^\tau
	\D \sigma \tW_\tau
	+ \theta_\gamma{}^\sigma \theta_\delta{}^\tau
	\D \sigma \theta_\tau{}^\mu
	\tW_\mu
\right)
\, dt \, dA_0,
\\
& = &
\int_\Omega
\int_{-h_0}^{h_0}
\left(
	\Lambda^{\alpha\beta\gamma\delta} e_{\alpha\beta} 
	B_{\gamma\delta}
	- t \, 
	\D \tau \D \sigma
	\left(
	\Lambda^{\alpha\beta\gamma\delta} e_{\alpha\beta} 
	\theta_\gamma{}^\sigma \theta_\delta{}^\tau
	\right)
\right)
\, \tomega
\, dt \, dA_0,
\\
& & \mbox{} +
\int_\Omega
\int_{-h_0}^{h_0}
\left(
	\Lambda^{\alpha\beta\gamma\delta} e_{\alpha\beta} 
	\theta_\gamma{}^\sigma \theta_\delta{}^\tau
	\D \sigma \theta_\tau{}^\mu
	-
	\D \sigma
	\left(
	\Lambda^{\alpha\beta\gamma\delta} e_{\alpha\beta} 
	\theta_\gamma{}^\sigma \ggg_\delta{}^\mu
	\right)
\right)
\tW_\mu
\, dt \, dA_0,
\end{eqnarray*}
Collecting the terms we obtain the equations
for the normal and the tangential components of the force
\begin{eqnarray*}
f^3
& = &
\int_{-h_0}^{h_0}
\left(
	\Lambda^{\alpha\beta\gamma\delta} e_{\alpha\beta} 
	B_{\gamma\delta}
	- t \, 
	\D \tau \D \sigma
	\left(
	\Lambda^{\alpha\beta\gamma\delta} e_{\alpha\beta} 
	\theta_\gamma{}^\sigma \theta_\delta{}^\tau
	\right)
\right)
\, dt
\\
f^\mu
& = & 
\int_{-h_0}^{h_0}
\left(
	\Lambda^{\alpha\beta\gamma\delta} e_{\alpha\beta} 
	\theta_\gamma{}^\sigma \theta_\delta{}^\tau
	\D \sigma \theta_\tau{}^\mu
	-
	\D \sigma
	\left(
	\Lambda^{\alpha\beta\gamma\delta} e_{\alpha\beta} 
	\theta_\gamma{}^\sigma \ggg_\delta{}^\mu
	\right)
\right)
\, dt
\end{eqnarray*}
Using (\ref {eqn:straincomponents}) again we can finally bring these
equations into the following form
%
%-------------------------------------------
% \input {force_eqns}
%%
%%
\begin {eqnarray}
\label{eq:Force1}
f^3 & = & A \, \omega
	+ \D \sigma \D \tau (\Abar^{\sigma\tau\mu\nu} \DDomega \mu \nu)
	- \D \sigma \D \tau (\Abarbar^{\sigma\tau} \omega)
	- \Abarbar^{\sigma\tau} \DDomega \sigma \tau
					 									% \nonumber 
\\
	& & \mbox{} 
	+ \Phi^\nu W_\nu
	+ \Phibar^{\mu\nu} \D \mu W_\nu
	- \D \sigma \D \tau (\Psi^{\mu\sigma\tau} W_\mu)
	- \D \mu \D \nu (\Psibar^{\sigma\tau\mu\nu} \D \sigma W_\tau)
					 									\nonumber 
\\% [-1.5ex]
\label{eq:Force2}
% \label{eq:Force}\\[-1.5ex]
f^\mu & = & \Omega^{\mu\nu}  \, W_\nu 
	+ \Omegabar^{\sigma\tau\mu} \, \D \sigma W_\tau
	- \D \sigma (\Omegabar^{\sigma\mu\tau} W_\tau) 
	- \D \nu (\Omegabarbar^{\sigma\tau\nu\mu} \D \sigma W_\tau)
														% \nonumber 
\\
	& & \mbox{} 
	+ \Phi^\mu \, \omega
	- \D \nu (\Phibar^{\nu\mu} \, \omega) 
	- \Psi^{\mu\sigma\tau} \DDomega \sigma \tau
	+ \D \nu (\Psibar^{\nu\mu\sigma\tau} \DDomega \sigma \tau)
					 									\nonumber 
\end {eqnarray}
%%
%%
%-------------------------------------------
%
where the coefficients are defined as follows:
%
%-------------------------------------------
% \input {force_coeffs}
%
\begin {eqnarray}
\label {eqnarray:coefficientsFirst}
A & = & \int \Lambda^{\alpha\beta\gamma\delta}
	% \theta_\alpha{}^\sigma b_{\sigma\beta}
	% \theta_\gamma{}^\tau b_{\tau\delta}
	B_{\alpha\beta} B_{\gamma\delta}
\, dt \\
\Abar^{\sigma\tau\mu\nu} & = &
	\int \Lambda^{\alpha\beta\gamma\delta}
	\theta_\alpha{}^\sigma \theta_\beta{}^\tau
	\theta_\gamma{}^\mu \theta_\delta{}^\nu
\, t^2 \, dt \\
\Abarbar^{\mu\nu} & = &
	\int \Lambda^{\alpha\beta\gamma\delta}
	% \theta_\alpha{}^\sigma b_{\sigma\beta}
	B_{\alpha\beta}
	\theta_\gamma{}^\mu \theta_\delta{}^\nu
\, t \, dt \\
\Phi^\mu & = &
	\int \Lambda^{\alpha\beta\gamma\delta}
	% \theta_\alpha{}^\sigma  b_{\sigma\beta}
	B_{\alpha\beta}
	\theta_\gamma{}^\tau \theta_\delta{}^\rho
	(\D \tau \theta_\rho{}^\mu)
\, dt \\
& = &                                   \nonumber
	\Abarbar^{\tau\rho}
	\D \tau b_\rho{}^\mu
\\
\Phibar^{\mu\nu} & = &
	\int \Lambda^{\alpha\beta\gamma\delta}
	% \theta_\alpha{}^\sigma b_{\sigma\beta}
	B_{\alpha\beta}
	\theta_\gamma{}^\mu \theta_\delta{}^\tau
	\theta_\tau{}^\nu
\, dt \\
\Psi^{\rho\mu\nu} & = &
	\int \Lambda^{\alpha\beta\gamma\delta}
	\theta_\alpha{}^\sigma \theta_\beta{}^\tau
	(\D \sigma \theta_\tau{}^\rho)
	\theta_\gamma{}^\mu \theta_\delta{}^\nu
\, t \, dt \\
& = &                                   \nonumber
	\Abar^{\sigma\tau\mu\nu}
	\D \sigma b_\tau{}^\rho
\\
\Psibar^{\sigma\tau\mu\nu} & = &
	\int \Lambda^{\alpha\beta\gamma\delta}
	\theta_\alpha{}^\sigma 
	% \theta_\beta{}^\lambda \theta_\lambda{}^\tau 
	\ggg_\beta{}^\tau
	\theta_\gamma{}^\mu
	\theta_\delta{}^\nu
\, t \, dt
\\
\Omega^{\mu\nu} & = &
	\int \Lambda^{\alpha\beta\gamma\delta}
	\theta_\alpha{}^\sigma \theta_\beta{}^\tau
	(\D \sigma \theta_\tau{}^\mu)
	\theta_\gamma{}^\lambda \theta_\delta{}^\rho
	(\D \lambda \theta_\rho{}^\nu)
\, dt
\\
& = &                                   \nonumber
	\Abar^{\sigma\tau\lambda\rho}
	\D \sigma b_\tau{}^\mu
	\D \lambda b_\rho{}^\nu
\\
\Omegabar^{\mu\nu\rho} & = &
	\int \Lambda^{\alpha\beta\gamma\delta}
	\theta_\alpha{}^\mu 
	% \theta_\beta{}^\sigma \theta_\sigma{}^\nu 
	\ggg_\beta{}^\nu
	\theta_\gamma{}^\tau
	\theta_\delta{}^\lambda
	(\D \tau \theta_\lambda{}^\rho)
\, dt
\\
& = &                                   \nonumber
	\Psibar^{\sigma\tau\mu\nu}
	\D \nu b_\mu{}^\rho
\\
\label {eqnarray:coefficientsLast}
\Omegabarbar^{\sigma\tau\mu\nu} & = &
	\int \Lambda^{\alpha\beta\gamma\delta}
	\theta_\alpha{}^\sigma 
	% \theta_\beta{}^\lambda \theta_\lambda{}^\tau 
	\ggg_\beta{}^\tau
	\theta_\gamma{}^\mu
	% \theta_\delta{}^\rho \theta_\rho{}^\nu
	\ggg_\delta{}^\nu
\, dt
\end {eqnarray}

%-------------------------------------------
%
%
\subsection{The Plate Equation}
For a flat shell surface $S$ the only non-zero coefficients are
$\Abar$ and $\Omegabbar$ ($\Psibar$ is zero because it is an integral
of an odd function on a symmetric interval),
and
\begin{eqnarray*}
\Lambda^{\sigma\tau\mu\nu}
& = &
	\frac {\lambda\mu} {\lambda + 2 \mu}
	h^{\sigma\tau} h^{\mu\nu}
	+ \mu
	h^{\sigma\mu} h^{\tau\nu}
\\
\Abar^{\sigma\tau\mu\nu} 
& = &
\frac 2 3 h_0^3
\,
\Lambda^{\sigma\tau\mu\nu}
\\
\Omegabbar^{\sigma\tau\mu\nu} 
& = &
2 h_0
\,
\Lambda^{\sigma\tau\mu\nu}
\end{eqnarray*}
The equations (\ref{eq:Force1}) -- (\ref{eq:Force2})
reduce to
%%  \begin {eqnarray*}
%%  f^3 & = & \D \sigma \D \tau (\Abar^{\sigma\tau\mu\nu} \DDomega \mu \nu)
%%  \\
%%  f^\mu & = & - \D \nu (\Omegabarbar^{\sigma\tau\nu\mu} \D \sigma W_\tau)
%%  \end {eqnarray*}
%
\begin {eqnarray*}
f_{normal} & = & \frac 2 3 h_0^3 D \, \Delta^2 \omega
\\
f_{tangential} & = & 2 h_0 D \, \D {}\div W
\end {eqnarray*}
where
$$
D = 2 \mu \frac {\lambda + \mu} {\lambda + 2 \mu}
\, .
$$
The first equation is the classical equation 
of an elastic plate (see \cite{Leissa1993}).
The second equation does not appear in classical plate analysis
where it is assumed that the plate undergoes only vertical motion
and the shearing forces are negligible.
%================================================================
% \input {ibm}

% \section {The Immersed Boundary Method}
\section {The Immersed Boundary Method}
\label {ChapterIBM}

%================================================================
% \input {NSeqns_symbols}

%some symbols used in the NS eqns

\def\boldu{{\bf u}}
\def\boldF{{\bf F}}
\def\boldf{{\bf f}}
\def\boldx{{\bf x}}
\def\boldX{{\bf X}}
\def\boldq{{\bf q}}

\def\dubydt{\frac {\partial \boldu} {\partial t}}
\def\dXbydt{\frac {\partial \boldX} {\partial t}}

%================================================================

The immersed boundary method of Peskin and McQueen is a general
numerical method for modeling an elastic material
immersed in a viscous incompressible fluid. 
Typically the immersed material has been modeled
as a collection of fibers. 
For details and
references to many applications see \cite {Peskin1994} 
and \cite{Peskin2002}.
This section outlines the immersed boundary method 
with the elastic material being modeled as a shell.

\subsection {The Equations of the Model}

For the description of the fluid we adopt a standard cartesian coordinate
system on $\Rthree$.
The immersed material is described in a different curvilinear 
coordinate system.
% In the immersed boundary method 
% the fluid is described in the standard cartesian coordinate system
% on $\Rthree$
% and the immersed material
% by a different curvilinear coordinate system.
The equations can be partitioned into three groups: the Navier Stokes
equations of a viscous incompressible fluid, the equations describing
the elastic material and the interaction equations.
Accordingly there are two distinct computational grids: one for the fluid
and the other for the immersed material. 
The purpose of the interaction equations is to communicate 
between these two grids.

We first turn to the Navier-Stokes equations.
Let $\rho$ and $\mu$ denote the density and the viscosity of
the fluid,
and let ${\boldu}({\boldx}, t)$ 
and $p({\boldx}, t)$ denote its velocity and pressure,
respectively.
% The viscous incompressible fluid is described by 
% the Navier-Stokes equations:
Then
\begin{eqnarray}
\label {eq: NS1}
\rho \left( \dubydt + \boldu \cdot \nabla \boldu \right)
        & = & \mbox{} - \nabla p + \mu \nabla^2 \boldu + \boldF
\\
\label {eq: NS2}
\nabla \cdot \boldu & = & 0
\end{eqnarray}
where $\boldF$ denotes the density of the body force 
acting on the fluid.
% the immersed material applies on the fluid. 
For example, if the immersed material
is modeled as a thin shell,
then $\boldF$ is a singular vector field,
which is zero everywhere, except possibly on the surface
representing the shell. 
 
Let $\boldX(\boldq, t)$ denote the position of the immersed material
in $\Rthree$.
% In the case of 
For a shell, $\boldq$ takes values in a domain 
$\Omega \subset \Rtwo$,
and $\boldX(\boldq, t)$ is a 1-parameter family of
surfaces indexed by $t$, 
i.e., $\boldX(\boldq, t)$ is the middle surface
of the shell at time $t$.
Let $\boldf(\boldq, t)$ denote the force density that 
the immersed material applies on the fluid.
% Then Newton's first law can be expressed as
Then
\begin{eqnarray}
\label {eq: fluidF}
\boldF(\boldx, t) = \int \boldf(\boldq, t)\delta(\boldx 
- \boldX(\boldq, t)) \,d\boldq,
\end{eqnarray}
where $\delta$ is the Dirac delta function on $\Rthree$.
This equation merely says that the fluid feels the force that the
material exerts on it, but it is important in the numerical method
where it is one of the interaction equations mentioned above.
The other interaction equation is the 
no slip condition for a viscous fluid:
\begin{eqnarray}
\label {eq: fluidU}
\dXbydt & = & \boldu (\boldX (\boldq, t), t) 
\nonumber \\
& = & \int \boldu(\boldx, t) \delta(\boldx - \boldX(\boldq, t))
\, d\boldx
\end{eqnarray}
\def\boldN{{\bf N}}
\def\boldW{{\bf W}}
%
%
% Equations (\ref{eq: fluidF}) -- (\ref{eq: fluidU})
% are the interaction equations mentioned above.
%
% Until now our description of the immersed boundary method did not depend
% on the nature of the material being modeled.
% To complete the description of the system we supplement
% the equations 
% (\ref{eq: NS1}) -- (\ref{eq: fluidU})
% with equations determining the material force $\boldf(\boldq, t)$.
% We will use the theory developed in chapter ??.
%
% position $\boldX(\boldq, t)$ which were derived 
%

\def\boldT{{\bf T}}
We will now complete the description of the system 
by writing down the third group of equations,
the equations describing the shell.
Let $\boldX_0(\boldq) $ denote the initial position,
which we take as our equilibrium reference configuration,
of the middle surface of the shell,
and let 
$\boldT_\mu = \boldT_\mu(\boldq) = \frac \partial {\partial q_\mu} $ 
($\mu = 1, 2$) be its tangent coordinate vector fields.
The displacement from the stationary configuration 
can be described in terms of the function $\omega(\boldq, t)$
and the vector field $\boldW(\boldq, t)$ defined by the following equation:
$$
	\boldX(\boldq, t) - \boldX_0(\boldq) = 
		\omega(\boldq, t) \boldN(\boldq) + \boldW(\boldq, t),
$$
where $ \boldN = \boldN(\boldq)$ denotes the normal to the surface 
$\boldX_0(\boldq) $
and $\boldW = W^\mu \boldT_\mu $ is tangent to $\boldX_0$.
We decompose the force $\boldf$ into its normal and tangential
components as well:
$$
\boldf = f^3 \boldN + f^\mu \boldT_\mu
$$
The components $f^i$ are related to $\omega$ and $W^\mu$ by the
shell equations 
(\ref {eq:Force1}) -- (\ref {eq:Force2}).
This completes the description of the fluid-shell system.

%
%-------------------------------------------------------------
%
\subsection {The Numerical Method}
\def\deltat{{\Delta t}}

Let $\deltat$ denote the duration of a time step.
It will be convenient to denote the time step by the superscript.
For example
$ \boldu^n(\boldx) = \boldu(\boldx, n \deltat). $
At the beginning of the $n$-th time step $\boldX^n$ and $\boldu^n$
are known.
Each time step proceeds as follows. 
\begin {enumerate}
\item
Compute the force $\boldf^n$ that the shell applies to the fluid.
\item
% Use the interpolation equations 
Use (\ref {eq: fluidF})
to compute the external force on the 
fluid $\boldF^n$.
\item
Compute the new fluid velocity $\boldu^{n + 1}$ from the Navier Stokes
equations.
\item
% Using the interpolation equations 
Use (\ref {eq: fluidU})
to compute the new position 
$\boldX^{n + 1}$ of the immersed material.
\end {enumerate}
The computation of the force in step 1 is explained in detail
in the next section.
We shall now describe in detail the computations in steps 2 --- 4,
beginning with the Navier-Stokes equations.

The fluid equations are discretized on a rectangular lattice of mesh
width $h$. 
We will make use of the following difference operators which act on
functions defined on this lattice:
\def\Dpi{{D^+_i}}
\def\Dmi{{D^-_i}}
\def\Dci{{D^0_i}}
\def\boldD0{{\bf D^0}}
\def\bolde{{\bf e}}
\begin{eqnarray}
\label {eq:differenceOpsFirst}
\Dpi \phi (\boldx) & = &
	\frac {\phi(\boldx + h \bolde_i) - \phi(\boldx)} h
\\
\Dmi \phi (\boldx) & = &
	\frac {\phi(\boldx) - \phi(\boldx - h \bolde_i)} h
\\
\Dci \phi (\boldx) & = &
	\frac {\phi(\boldx + h \bolde_i) - \phi(\boldx - h \bolde_i)} {2h}
\\
\label {eq:differenceOpsLast}
\boldD0 & = &
	(D^0_1, D^0_2, D^0_3)
\end{eqnarray}
where $ i = 1, 2, 3 $, and $\bolde_1$, $\bolde_2$, $\bolde_3$ form
an orthonormal basis of $\Rthree$.

In step 3 we use the already known $\boldu^n$ and  $\boldF^n$ to
compute  $\boldu^{n + 1}$ and $p^{n+1}$ by solving the following
linear system of equations:

\begin {eqnarray}
\label {eq: discreteNS1}
\rho \left( \frac {\boldu^{n + 1} - \boldu^n} {\Delta t}
        + \sum_{k = 1}^3 u_k^n D_k^\pm\boldu^n \right)
        & = &
\mbox{} 
        - {\bf D}^0 p^{n + 1}
        + \mu \sum_{k = 1}^3 D^+_k D^-_k \boldu^{n + 1}
        + \boldF^n
\\
\label {eq: discreteNS2}
{\bf D}^0 \cdot \boldu^{n + 1} & = & 0
\end {eqnarray}
Here 
% in equation (\ref{eq: discreteNS1}) 
% the notation 
$u_k^n D_k^\pm$
stands for upwind differencing:
$$
u_k^n D_k^\pm = \left\{
\begin {array} {ll}
u_k^n D_k^- 
& 
u_k^n > 0
\\
u_k^n D_k^+ 
&
u_k^n < 0
\end {array}
\right.
$$
Equations (\ref{eq: discreteNS1}) -- (\ref{eq: discreteNS2})
are linear constant coefficient difference equations 
and, therefore, can be solved efficiently with the use of the Fast
Fourier Transform algorithm.

We now turn to the discretization of equations
(\ref {eq: fluidF}), (\ref {eq: fluidU}).
Let us assume, for simplicity, that $\Omega \subset \Rtwo$
is a rectangular domain over which all of the quantities related
to the shell are defined. We will assume that this domain is
discretized with mesh widths $\Delta q_1$, $\Delta q_2$
and the computational lattice for $\Omega$ is the set
\def\boldQ{{\bf Q}}
$$
\boldQ = \left\{(i_1\Delta q_1, i_2 \Delta q_2) \ | 
\ i_1 = 1 \ldots n_1,
\ \ 
i_2 = 1 \ldots n_2 \right\}.
$$

In step 2 the force $\boldF^n$ is computed using the following equation.
\begin {eqnarray}
\label {eq:discreteF}
\boldF^n(\boldx) = \sum_{\boldq \in \boldQ}
	\boldf^n(\boldq)\delta_h(\boldx 
	- \boldX^n(\boldq)) \Delta\boldq
\end {eqnarray}
where 
% $\boldq = (q_1, q_2)$,
$ \Delta\boldq = \Delta q_1 \Delta q_2 $
and 
$\delta_h$ is a smoothed approximation to the Dirac delta function
on $\Rthree$ described below. 
% Here the summation is over the values
% $ \boldq = (i \Delta q_1, j \Delta q_2) $ 
% where $i$ and $j$ are integers.

Similarly, in step 4 updating the position of the immersed material 
$\boldX^{n + 1}$ is done using the equation
\begin {eqnarray}
\label {eq:discreteX}
\boldX^{n + 1}(\boldq) = \boldX^n(\boldq) + \Delta t 
\sum_\boldx \boldu^{n + 1}(\boldx) \delta_h(\boldx - \boldX^n(\boldq)) h^3
\, ,
\end {eqnarray}
where the summation is over the lattice
$ \boldx = (hi, hj, hk) $, where $i$, $j$ and $k$ are integers.

The function $\delta_h$ which is used in 
(\ref {eq:discreteF}) and (\ref {eq:discreteX}),
is defined as follows:
$$
\delta_h(\boldx) = 
	h^{-3}\phi(\frac {x_1} h) \phi(\frac {x_2} h) \phi(\frac {x_3} h)
\, ,
$$
where
$$
\phi(r) = \left\{
\begin {array} {cl}
\frac 1 8 (3 - 2 |r| + \sqrt{1 + 4|r| - 4 r^2}) 
&
|r| \le 1
\\
\frac 1 2 - \phi(2 - |r|)
&
1 \le |r| \le 2
\\
0 & 2 \le |r|
\end {array}
\right.
$$
For an explanation of the construction of $\delta_h$ see 
\cite {Peskin1994}.

%
%---------------------------------------------------------------
%
\subsection {Computation of the Elastic Force}

The force $\boldf^n$, that the shell applies to the fluid is computed
by discretizing equations (\ref {eq:Force1}) -- (\ref {eq:Force2}).
We shall now describe how to discretize these equations and how
to compute various geometric quantities,
such as the Christoffel symbols,
the second fundamental form, etc. 
These basic geometric quantities,
as well as the coefficients 
(\ref {eqnarray:coefficientsFirst}) 
-- (\ref {eqnarray:coefficientsLast})
that depend on them,  can
be computed once in the initialization step of the algorithm and
stored for subsequent use.

% Before we turn to the computation of the force $\boldf^n$ in step 1
% some preliminaries.
% we describe the preliminary computations of various geometric
% quantities.
A covariant derivative of a $(p,q)$-tensor $A$ is defined by
\begin {equation}
\label {eq:TensorDerivative}
\D \alpha A_{\mu_1 \ldots \mu_p}^{\nu_1 \ldots \nu_q} =
\partial_\alpha A_{\mu_1 \ldots \mu_p}^{\nu_1 \ldots \nu_q}
+ \sum_{k = 1}^q \Gamma_{\alpha\sigma}^{\nu_k}
	A_{\mu_1 \ldots \mu_p}^{\nu_1 \ldots \sigma \ldots \nu_q}
- \sum_{m = 1}^p \Gamma_{\alpha\mu_m}^\sigma
	A_{\mu_1 \ldots \sigma \ldots \mu_p}^{\nu_1 \ldots \nu_q}
\, ,
\end {equation}
where
$$
\Gamma^\lambda_{\mu \nu} = \frac 1 2 g^{\lambda \sigma}
\left( g_{\mu \sigma,\nu} + g_{\sigma \nu, \mu} - g_{\mu \nu, \sigma}
\right)
$$
are the Christoffel symbols,
$g_{\mu\nu} = \boldT_\mu \cdot \boldT_\nu $
is the metric of the surface $\boldX_0$,
$g^{\mu\nu}$ is the inverse of $g_{\mu\nu}$,
and comma denotes partial differentiation
(i.e., $\phi_{,\alpha} = \partial_\alpha \phi$).
% We discretize the computation of a covariant derivative by
To simplify the discretization of the force equations
we introduce a single difference operator $D_\alpha$
which uses center differencing in the interior of the domain and
either forward or backward differencing on the boundary:
$$
D_\alpha \phi (\boldq) =
\left\{
\begin {array} {ll}
D^+_\alpha \phi (\boldq)
&
q_\alpha = \Delta q_\alpha
\\
D^0_\alpha \phi (\boldq)
&
2 \Delta q_\alpha \le q_\alpha \le (n_\alpha - 1) \Delta q_\alpha
\\
D^-_\alpha \phi (\boldq)
&
q_\alpha = n_\alpha \Delta q_\alpha
\end {array}
\right.
$$
Here the operators $ D^+_\alpha $, $D^0_\alpha$ and $D^-_\alpha$
are defined on the lattice $\boldQ$ analogously to the definitions
(\ref {eq:differenceOpsFirst}) -- (\ref {eq:differenceOpsLast}).
We define the discrete covariant derivative of a $(p,q)$-tensor $A$
by
\begin {equation}
\DT \alpha A_{\mu_1 \ldots \mu_p}^{\nu_1 \ldots \nu_q} =
D_\alpha A_{\mu_1 \ldots \mu_p}^{\nu_1 \ldots \nu_q}
+ \sum_{k = 1}^q \Gamma_{\alpha\sigma}^{\nu_k}
	A_{\mu_1 \ldots \mu_p}^{\nu_1 \ldots \sigma \ldots \nu_q}
- \sum_{m = 1}^p \Gamma_{\alpha\mu_m}^\sigma
	A_{\mu_1 \ldots \sigma \ldots \mu_p}^{\nu_1 \ldots \nu_q}
\, .
\end {equation}
% 
% The discretization of the force equations is now straightforward:
% In practice we compute covariant derivatives by
% replacing the $\partial_\alpha$ operator with the operator
% $ D_\alpha $ defined as follows:

%
%
The following computations are carried out in the initialization stage
of the algorithm:
\begin {enumerate}
\item
Compute the tangent vector fields 
$
\boldT_\alpha = D_\alpha \boldX_0
$.
% $\alpha = 1, 2$.
\item
Compute the unit normal vector field 
$
\boldN = \frac {\boldT_1 \times \boldT_2}
					{\left| \boldT_1 \times \boldT_2 \right|}
\, .
$
\item
Compute the metric
$
g_{\mu\nu} = \boldT_\mu \cdot \boldT_\nu
$
and its  inverse $g^{\mu\nu}$.
\item
Compute the second fundamental form:
$
b_{\mu\nu} = D_\mu \boldN \cdot \boldT_\nu
$
.
\item
Compute the Christoffel symbols:
$$
	\Gamma^\lambda_{\mu\nu} = \frac 1 2 g^{\sigma\lambda}
		(D_\nu g_{\mu\sigma} + D_\mu g_{\sigma\nu}
			- D_\sigma g_{\mu\nu})
\, .
$$
\item
Compute the derivative of the second fundamental form:
$$
\DT \alpha b_\beta{}^\gamma = 
D_\alpha b_\beta{}^\gamma
	+ \Gamma_{\alpha\sigma}^\gamma b_\beta{}^\sigma
	- \Gamma_{\alpha\beta}^\sigma b_\sigma{}^\gamma
\, .
$$
\item
Compute the coefficients (\ref {eqnarray:coefficientsFirst}) --
(\ref {eqnarray:coefficientsLast}).
In practice instead of evaluating the integrals in these expressions
it is easier to
assume that the shell is sufficiently thin and to drop the 
terms involving high powers of the thickness.
% See the Appendix for details.
\end{enumerate}
%
% It would be too tedious to evaluate the integrals explicitly by hand.
% It could be more convenient to do so with the help of some software
% capable of performing symbolic integration.
% Alternatively one could perform numerical integration.
% However if we
% assume the shell is thin, we can neglect the powers of $t$ greater
% than 3
% and compute the integrals explicitly.
%
% Let $2h$ denote the thickness of the shell.
% Assume $C^{\alpha\beta\gamma\delta}$ is constant through the thickness.
% Integrating (...) and keeping 
% gives:
%
%
%

% Having computed and stored the coefficients during the initialization,
% we use a straightforward discretization 
% of the equations (\ref {eq:Force1}) - (\ref {eq:Force2})
% to compute the force at each time step:
%
The computation of the force during each time step proceeds as follows:
\begin {enumerate}
\item
Decompose the displacement into its tangential and normal components:
$$
\omega = (\boldX^n - \boldX_0) \cdot \boldN
$$
$$
W_\alpha = (\boldX^n - \boldX_0) \cdot \boldT_\alpha
$$
\item
Compute $D_\mu \omega$, the components of $d\omega$, the 1-form
dual to the gradient of $\omega$.
\item
Compute $\DT \mu D_\nu  \omega$, the components of the hessian 
of $\omega$.
\item
Compute the components $\DT \mu W_\nu$.
\item
Compute the normal and the tangential components of the force
using the following discretization of the equations
(\ref {eq:Force1}) - (\ref {eq:Force2}):
\begin {eqnarray}
f^3 & = & A \, \omega
	+ \DT \sigma \DT \tau (\Abar^{\sigma\tau\mu\nu} \DT\mu D_\nu \omega)
	- \DT \sigma \DT \tau (\Abarbar^{\sigma\tau} \omega)
	- \Abarbar^{\sigma\tau} \DT\sigma D_\tau \omega 
					 									\nonumber \\
	& & \mbox{} 
	+ \Phi^\nu W_\nu
	+ \Phibar^{\mu\nu} \DT \mu W_\nu
	- \DT \sigma \DT \tau (\Psi^{\mu\sigma\tau} W_\mu)
	- \DT \mu \DT \nu (\Psibar^{\sigma\tau\mu\nu} \DT \sigma W_\tau) \\
f^\mu & = & \Omega^{\mu\nu}  \, W_\nu 
	+ \Omegabar^{\sigma\tau\mu} \, \DT \sigma W_\tau
	- \DT \sigma (\Omegabar^{\sigma\mu\tau} W_\tau) 
	- \DT \nu (\Omegabarbar^{\sigma\tau\nu\mu} \DT \sigma W_\tau)
														\nonumber \\
	& & \mbox{} 
	+ \Phi^\mu \, \omega
	- \DT \nu (\Phibar^{\nu\mu} \, \omega) 
	- \Psi^{\mu\sigma\tau} \DT\sigma D_\tau \omega
	+ \DT \nu (\Psibar^{\nu\mu\sigma\tau} \DT\sigma D_\tau \omega)
\end {eqnarray}

\end {enumerate}

To complete the computation, we 
express the force in cartesian coordinates using
$$
\boldf = f^3 \boldN + f^\mu \boldT_\mu.
$$

%================================================================
% \input {application}

\section {An Application: Modeling a Shell Immersed in Fluid}
\label {ChapterApplication}

In this chapter we present a model shell which is a prototype
of a piece of the basilar membrane of the cochlea.
The immersed boundary method for shells was implemented using
the C programming language and numerical experiments were carried
out with this shell. 
We describe the model
that was used,
% in the numerical experiments,
examine the convergence of the algorithm in this case
and discuss the results.

\subsection {The Model Shell}

We construct a shell similar to a piece of the basilar membrane. 
Our shell will be approximately one seventh of the length
of the basilar membrane in the human cochlea
and it will have more curvature.
Introducing more curvature allows us to pack a longer shell into
a cube of fluid of a given size. 
% Because of computer limitations (see Chapter \ref {chapter:conclusion}
% for a detailed discussion),
Choosing a small fluid cube enables us
to achieve a better numerical resolution of the fluid.
% because at present it is not practical to perform computations
% on a fluid grid of more than $128^3$ points.
%
% \begin{figure}
% \setlength{\epsfxsize}{4.5in}
% % \vspace {0.2in}
% % \hspace{0.5in}
% \epsffile{../shell.ps}
% \caption{The model shell.}
% \label {fig:shell}
% \end{figure}
%
%
% We choose the width to grow linearly and the thickness of the shell,
% as well as its material properties,
% are chosen to make its compliance similar
% to that of a piece of basilar membrane.
%

The surface is a narrow helicoidal strip defined as follows.
Consider the curve
$$
\gamma(t) = 
	(R \cos(\alpha t), R \sin(\alpha t), H \alpha t) 
$$
where 
$R$, $\alpha$ and $H$ are constants.
They are specified, along with the other parameters that we use below,
in Table \ref {tab:parameters}.
%
%--------------------------------table--------------------------
\begin {table}
\caption {Parameters of the numerical test model}
\vspace {0.2in}
\begin {tabular} {|l|l|}
\hline
% \\
$a = 0.1$ cm
	&
	length of the side of the fluid cube
\\
$N = 32, 64, 128$
	&
	size of the fluid grid is $N^3$
\\
$h = \frac a N$
	&
	fluid mesh width
\\
$\rho = 1.034$ cm~g$^{-3}$
	&
	fluid density
\\
$\nu = 0.0197$ g~cm$^{-1}$~sec$^{-1}$
	&
	fluid viscosity
\\
$L_{BM} = 3.5$ cm
	& 
	length of the basilar membrane
\\
$L = 0.5$ cm
	&
	length of the model shell
\\
$w_0 = 0.015$ cm
	&
	width of the basilar membrane at the base
\\
$w_1 = 0.056$ cm
	&
	width of the basilar membrane at the apex
\\
$w(q_1) = w_0 + \frac {q_1} {L_{BM}} (w_1 - w_0)$
	&
	width of the model shell
\\
$h(q_1) \approx 0.001 (1.0 + 5.0~q_1) $ cm
	&
	thickness of the model shell
\\
$\alpha = 1.8 \pi / L$
	&
\\
$R = \frac 1 {30} $ cm
	&
\\
$ H = 0.01 $ cm
	&
\\
$n_1 = 1280 \frac N  {128}$
	&
	first dimension of the surface grid
\\
$n_2 = 48 \frac N  {128}$
	&
	second dimension of the surface grid
\\
$\Delta q_1 = \frac L {n_1 - 1}$
	&
	surface mesh width
\\
$\Delta q_2(q_1) = \frac {w(q_1)} {n_2 - 1}$
	&
	surface mesh width
\\
$\lambda = 26197503.0$ g~cm$^{-1}$~sec$^{-2}$
	&
	first Lam\'{e} coefficient
\\
$\mu = 523950.0$ g~cm$^{-1}$~sec$^{-2}$
	&
	second Lam\'{e} coefficient
\\
$\Delta t = 0.5,~1.0,~2.0,~4.0 \times 10^{-8} $ sec
	&
	time step
\\
$T_0 = 2.0 \times 10^{-6} $ sec
	&
	total simulation time
\\
\hline
\end {tabular}
\label {tab:parameters}
\end {table}
%--------------------------------table--------------------------
%
%
%
The vector field
$$
N(t) = (-\cos(\alpha t), -\sin(\alpha t), 0)
$$
is a unit normal field along the curve $\gamma$.
We define the surface of the shell by the following equation:
\begin {eqnarray*}
\boldX(q_1, q_2) & = & \gamma(q_1) 
	+ \left(q_2 - \frac {w(q_1)}  2\right) N(q_1),
\\
& &  0 \le q_1 \le L,
\qquad \frac {-w(q_1)} 2 \le q_2 \le \frac {w(q_1)} 2 
\end {eqnarray*}
%
%
%
% (see Figure \ref {fig:shell}).
The width $w$ is chosen to grow linearly, similar to the width
of the basilar membrane (see Table \ref {tab:parameters}).
This surface is discretized on a grid of size $n_1 \times n_2$
% using the following equations:
as follows:
\begin {eqnarray*}
\boldX_{k_1, k_2} & = & 
	% \gamma_{k_1} 
	\gamma(k_1 \Delta q_1)
	+ \left(k_2 \, \Delta q_2(k_1 \Delta q_1)
	- \frac {w(k_1 \Delta q_1)}  2
	\right)
	% N_{k_1}
	N(k_1 \Delta q_1)
% \\
% \gamma_{k_1} & = & \gamma(k_1 \Delta q_1)
% \\
% N_{k_1} & = & N(k_1 \Delta q_1)
\\
& &
	\qquad
	k_1 = 1, ..., n_1,
	\qquad
	k_2 = 1, ..., n_2.
\end {eqnarray*}

\begin{figure}
% \vspace {0.2in}
% \hspace{0.5in}
% \includegraphics[height=3.0in]{h.ps}
\setlength{\epsfxsize}{3.0in}
\epsffile{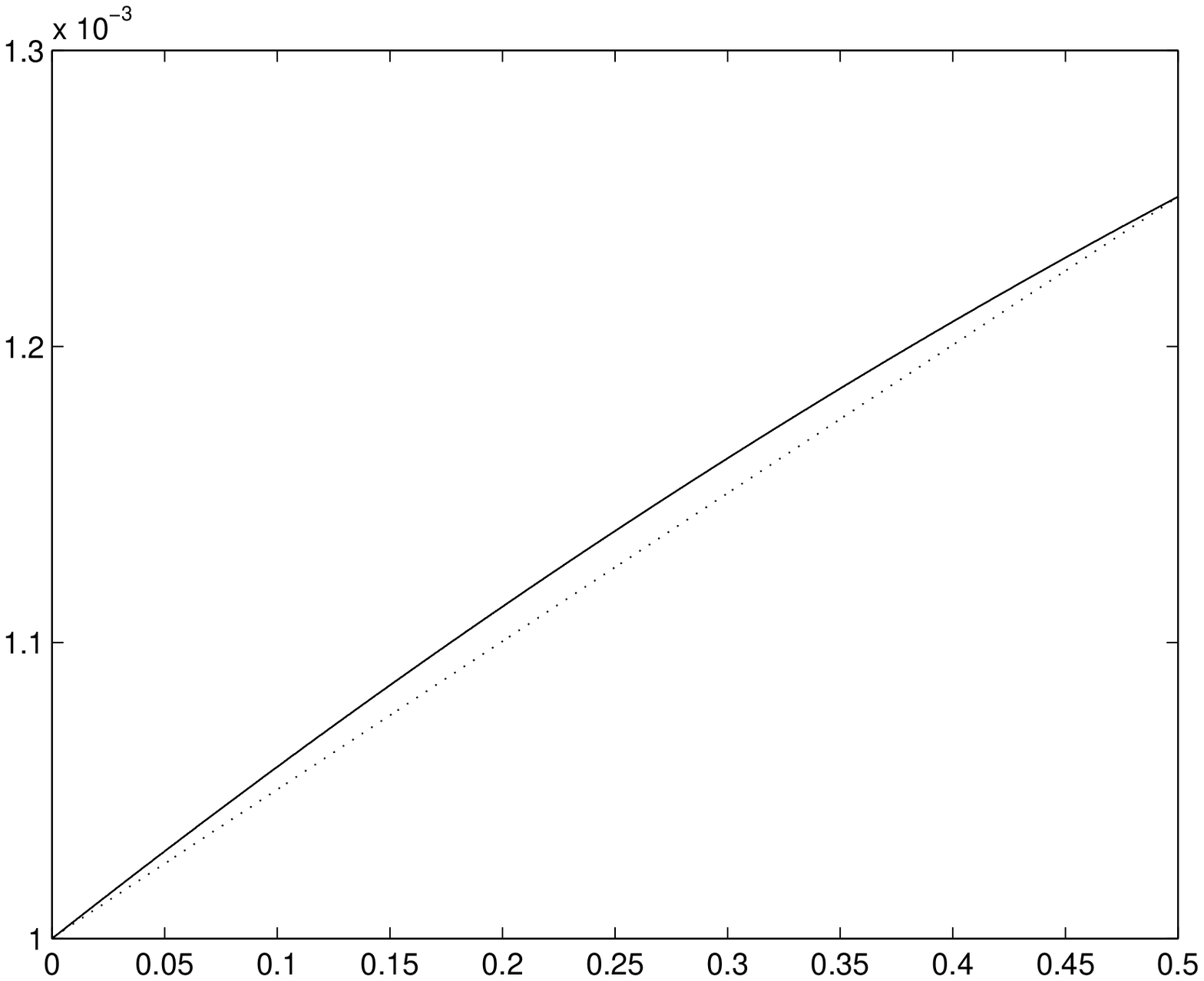}
\caption{The thickness of the model shell.}
\label {fig:thickness}
\end{figure}

The compliance of an elastic shell is defined as the amount
of volume displaced per unit pressure difference across the shell.
Von Bekesy found that the compliance per unit length of the basilar 
membrane varies exponentially along the basilar membrane
as $e^{c q_1}$, 
where $c^{-1} = 0.7$ cm
(see \cite {vonBekesy}).
To achieve this compliance in the model shell we use the shell
equations to estimate
the required thickness. 
(Notice that nothing changes in the derivation of the shell equations 
when the thickness $h_0$
is assumed to be a function on the middle surface,
rather than being a constant. 
The thickness enters the equations only in the definition of the
coefficients
(\ref {eqnarray:coefficientsFirst}) --
(\ref {eqnarray:coefficientsLast})).
We obtain an estimated thickness
%
% Consequently, we choose the  thickness of the shell to be:
$$
h(q_1) = 0.001 ~ \left(1 + \frac 2 3 q_1\right)^{\frac 5 3}
				10^{-\frac 2 9 q_1} \textrm { cm,}
$$
which on the interval $0 \le q_1 \le 0.5$ is very close to
a straight line (Figure \ref {fig:thickness}).
The above choices of the width, the thickness and the L\'{a}me
coefficients (see Table \ref {tab:parameters}) yield an
estimated compliance similar to the one measured by von Bekesy
for the basilar membrane.

%
%-----------------------------------------------------------------------
%
\subsection {The Numerical Experiments}

The immersed boundary method for shells was implemented
using the C programming language.
The program ran on the Cray C-90 at the Pittsburgh 
Supercomputing Center as well as on Silicon Graphics 
workstations.

The numerical experiments are set up as follows.
The shell $\boldX$ is placed in a periodic cube of fluid,
i.e., a three-dimensional torus.
The length of this cube's side is $a = 0.1$ cm,
so it is much smaller than the cochlea, whose volume is
approximately 1 cm$^3$.
The fluid density $\rho$, and the fluid viscosity $\nu$,
are chosen equal to the corresponding parameters of the cochlear
fluid.
The fluid equations are discretized on a cubic grid of $N^3$ points
and 
the grid size for the shell is chosen correspondingly
such that its 
mesh widths, $\Delta q_1$ and $\Delta q_2$, 
% would be slightly smaller than 
are approximately equal to half of the fluid mesh width $h$
(Table \ref {tab:parameters}).
This choice is necessary to prevent the fluid from leaking 
through the shell.
The edges of the shell's boundary are clamped
to a fixed location in space with two rows of springs along
each edge.
% The boundary of the shell is clamped with two rows of springs.
%
%
% In addition to varying the size of the computational grids,
% we have varied the time step $\Delta t$.
%

At time 0 the system is perturbed with a force impulse in the fluid
and the simulation is run for a fixed period of time 
$T_0 = 2.0 \times 10^{-6} $ sec.
The force impulse acting at time 0
is represented by a constant singular vector field 
in the vertical direction
defined on the horizontal plane $z = 0.1~\mbox{cm}$ 
(where $z$ is the vertical coordinate)
with the force density of
$4 \times 10^{-7}~\mbox{g cm}^{-1} \mbox{sec}^{-2}$.
Since the force source is broadly distributed through a
horizontal plane within the fluid, wave propagation within
the basilar membrane can only arise as a consequence of its material
properties
and cannot be the result of the location of the force source.

Immersed boundary computations typically require large scale
computing resources.
Because of time and storage limitations, it is not possible to
conduct an extensive empirical study of the algorithm's convergence.
At present it is not practical to implement
the method with a fluid grid of more than $128^3$ points.
The experiment has been repeated with 
$N = 32, 64 \textrm{ and } 128$ 
and time steps $\Delta t = 0.5,~1.0,~2.0
\textrm{ and } 4.0 \times 10^{-8}$ sec.
% The Courant-Friedrichs-Lewy condition forces a choice of such 
Numerical stability condition forces a choice of such 
small time steps.
On the other hand, reducing the time step further
may result in a significant machine precision error.
% Thus, twelve different sets of data were obtained.

%
%-----------------------------------------------------------------------
%
\subsection {Convergence of the Algorithm}

Let $\boldX_1(t)$ and $\boldX_2(t)$ denote the position of the shell
at time $t$ obtained in two different numerical experiments.
To measure the relative difference at time $t$
we calculate
\begin{equation}
E(t) = \frac {|\boldX_1(t) - \boldX_2(t)|_p} 
		{|\boldX_1(t) - \boldX_1(0)|_p}
\label {eq:error}
\end{equation}
where the norms are $L^p$-norms with $p = 1, 2 \textrm{ and } \infty$.
% as a measure of the relative difference between two tests at time $t$.
% We have used $L_1, L_2$ and max-norms in (\ref {eq:error}).

%
\begin{figure}[!htpb]
% \vspace{0.2in}
% \hspace{0.5in}
% \includegraphics[height=3.0in]{12_88.ps}
\setlength{\epsfxsize}{3.0in}
\epsffile{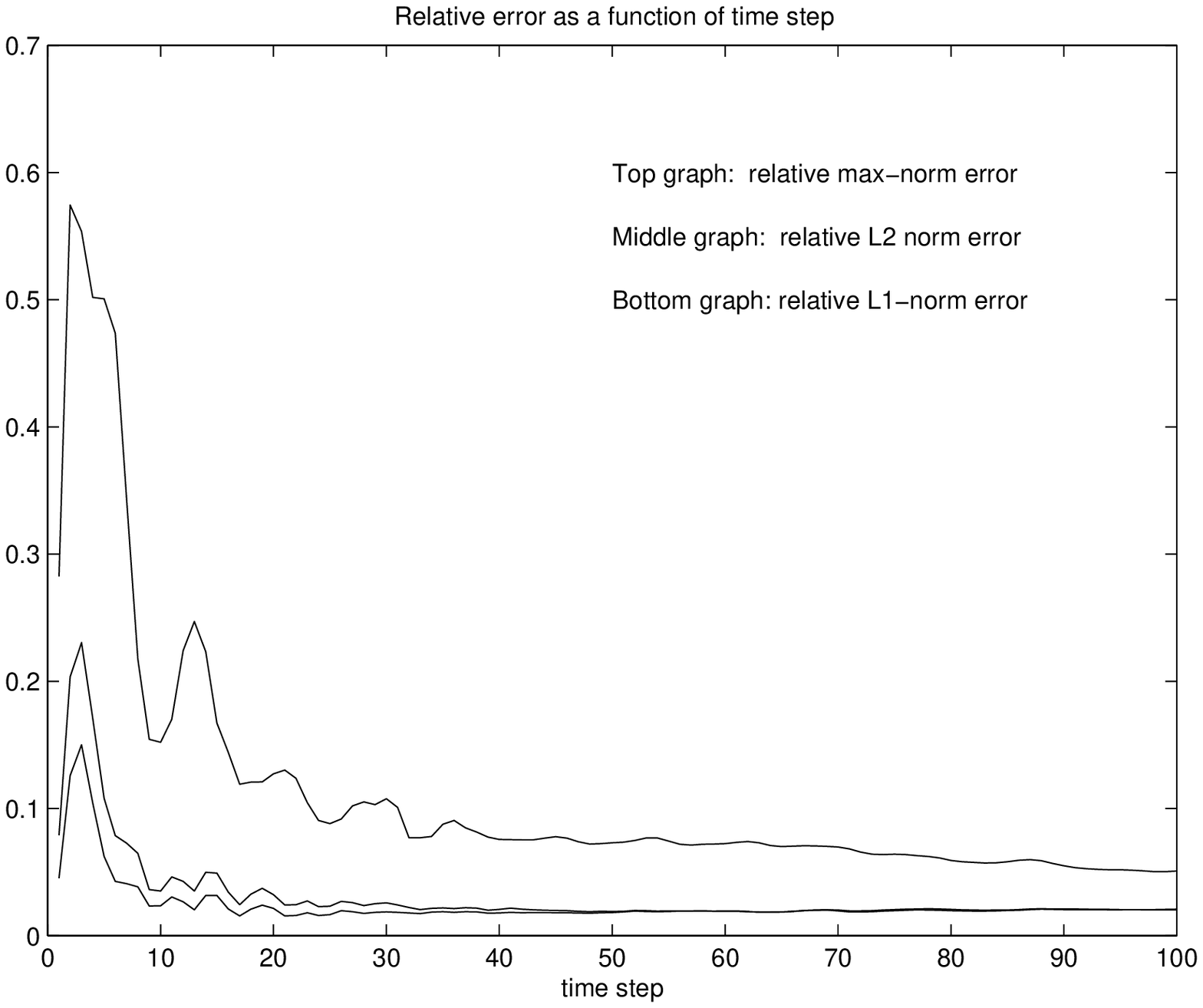}
\caption{The relative difference between $\boldX_{1,128}$
and $\boldX_{2,128}$.}
\label {errora}
\end{figure}

\begin{figure}[!htbp]
% \vspace{0.2in}
% \hspace{0.5in}
% \includegraphics[height=3.0in]{55_84.ps}
\setlength{\epsfxsize}{3.0in}
\epsffile{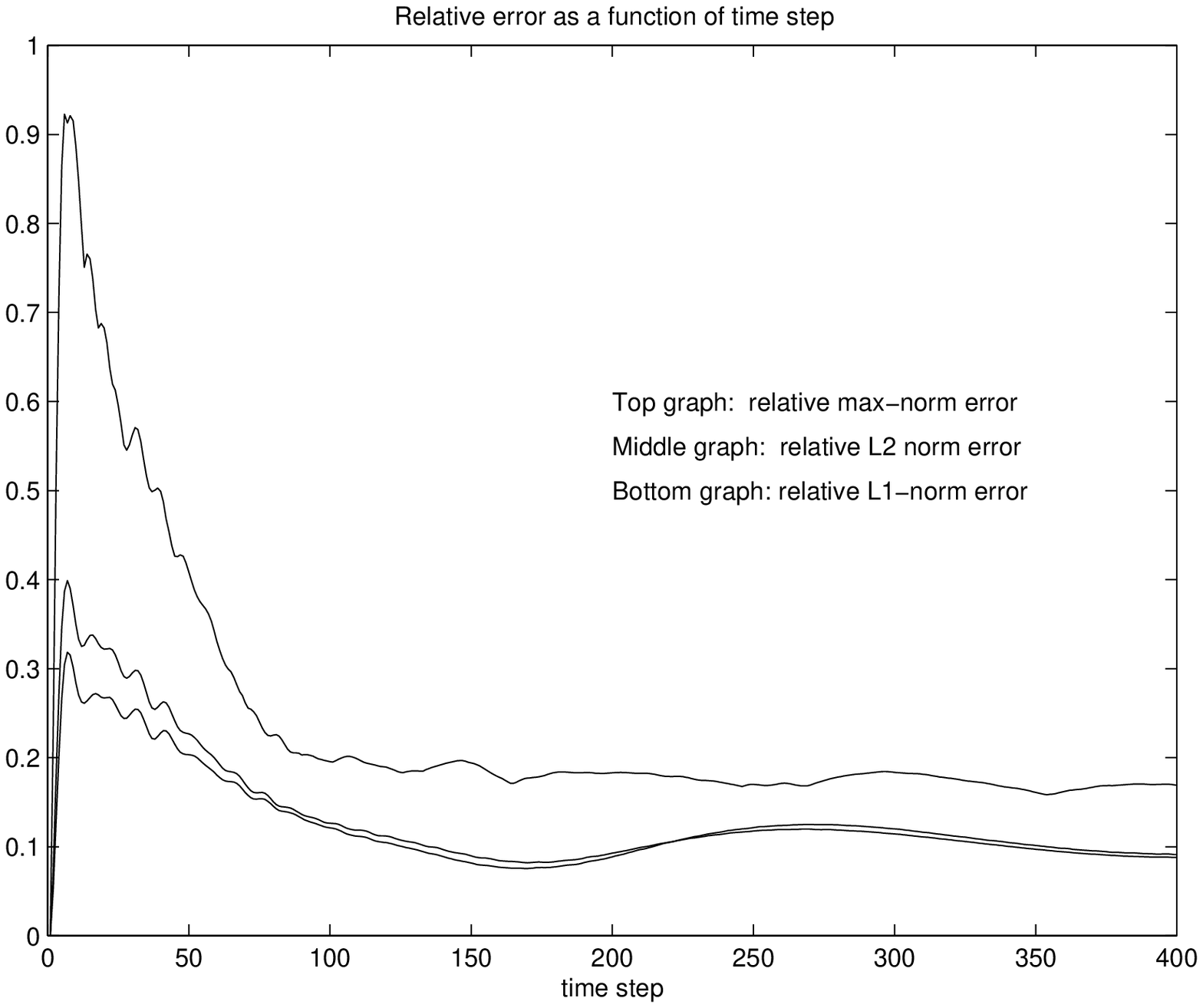}
\caption{The relative difference between $\boldX_{0.5,128}$
and $\boldX_{0.5,64}$.}
\label {errorb}
\end{figure}

For any two experiments
the function $E$ has been found to decrease initially very sharply,
stabilizing in the second half of the time interval.
% This was observed for any two data sets.
Two typical graphs are shown in
Figure \ref {errora}
and
Figure \ref {errorb}.
Therefore, we will measure the difference
between two computed solutions on the time interval 
$[\frac {T_0} 2,~ T_0]$ using the following
space-time norms:
$$
\|\boldX_1 - \boldX_2\|_p = \sum_{t \in [0.5~T_0, ~T_0]}
|\boldX_1(t) - \boldX_2(t)|_p
$$
where the summation is taken over the set of 25 time values
common to all of the tests performed
and the $L^p$-norms are computed on the ``common'' grid of
$257 \times 12$ points.

The computations were performed with double precision floating point
on a Silicon Graphics computer.
The results of the measurements
% Measurements' results 
are summarized in Table \ref {tab:NormComparison}
and Table \ref {tab:convergence}.
These results indicate that the algorithm converges
when the time step $\Delta t$ is linearly proportional
to the mesh width $h$ and both tend to zero.
\begin {table}[!htbp]
\caption {Norm comparison between different numerical tests.}
\vspace {0.2in}
\begin {tabular}{|rcr|l|l|l|}
\hline
\multicolumn {3} {|c|} {Tests compared}
&
$L^1$-norm
&
$L^2$-norm
&
$L^\infty$-norm
\\
\hline
1.0/128 & -- & 2.0/64
&
$ 9.6290 \times 10^{-08} $
&
$ 1.8547 \times 10^{-09} $
&
$ 9.7749 \times 10^{-11} $
\\
2.0/64 & -- & 4.0/32
&
$ 2.3984 \times 10^{-07} $
&
$ 4.6353 \times 10^{-09} $
&
$ 1.8326 \times 10^{-10} $
\\
0.5/128 & -- & 1.0/64
&
$ 9.2989 \times 10^{-08} $
&
$ 1.8086 \times 10^{-09} $
&
$ 9.5158 \times 10^{-11} $
\\
1.0/64 & -- & 2.0/32
&
$ 2.4082 \times 10^{-07} $
&
$ 4.7011 \times 10^{-09} $
&
$ 1.9287 \times 10^{-10} $
\\
\hline
\end {tabular}
\label {tab:NormComparison}
\end {table}
Let $\boldX_{\delta, \overline {N}}$ denote the position of the shell
computed with $ \Delta t = \delta \times 10^{-8} $ and $N = \overline {N}$.
Using the data in Table \ref {tab:NormComparison}
we have two estimates of the order of convergence
of the algorithm:
$$
r_1 = \log_2 \left( 
		\frac {\| \boldX_{1,128} - \boldX_{2,64} \|_p} 
		{\| \boldX_{2,64} - \boldX_{4,32} \|_p }
		\right)
$$
and
$$
r_2 = \log_2 \left(
		\frac {\| \boldX_{0.5,128} - \boldX_{1,64} \|_p} 
		{\| \boldX_{1,64} - \boldX_{2,32} \|_p }
		\right)
$$
The values of $r_1$ and $r_2$ are given in
Table \ref {tab:convergence}.
As the time step becomes smaller, the machine precision error
becomes more significant leading to a slower rate of convergence.
This is apparently the reason that we have $r_2 < r_1$.
\begin {table}
\caption {Convergence rate estimates.}
\vspace {0.2in}
\begin {tabular}{|l|l|l|l|}
\hline
&
$L^1$-norm
&
$L^2$-norm
&
$L^\infty$-norm
\\
\hline
$r_1$
&
1.3729          & 1.3781          & 1.0192
\\
$r_2$
&
1.3166          & 1.3215          & 0.9068
\\
\hline
\end {tabular}
\label {tab:convergence}
\end {table}

%
%
%
%----------------------------------------------------------
%
\subsection {The Traveling Wave}

The test model described above was run with $N = 128$
and time step $\Delta t = 1.0 \times 10^{-8}$ seconds
for a total of more than 2400 time steps.
Although this  model is still far from a realistic model
of the cochlea,
several qualitative features, characteristic of the cochlear
wave mechanics, were already observed in this experiment.

\begin{figure}[h]
\begin{center}
$\begin{array}{c@{\hspace{0.1in}}c}
% \multicolumn{1}{l}{\mbox{\bf (a)}} &
% \multicolumn{1}{l}{\mbox{\bf (b)}} \\ [-0.53cm]
\vspace{-3cm}
\epsfxsize=2.5in
\epsffile{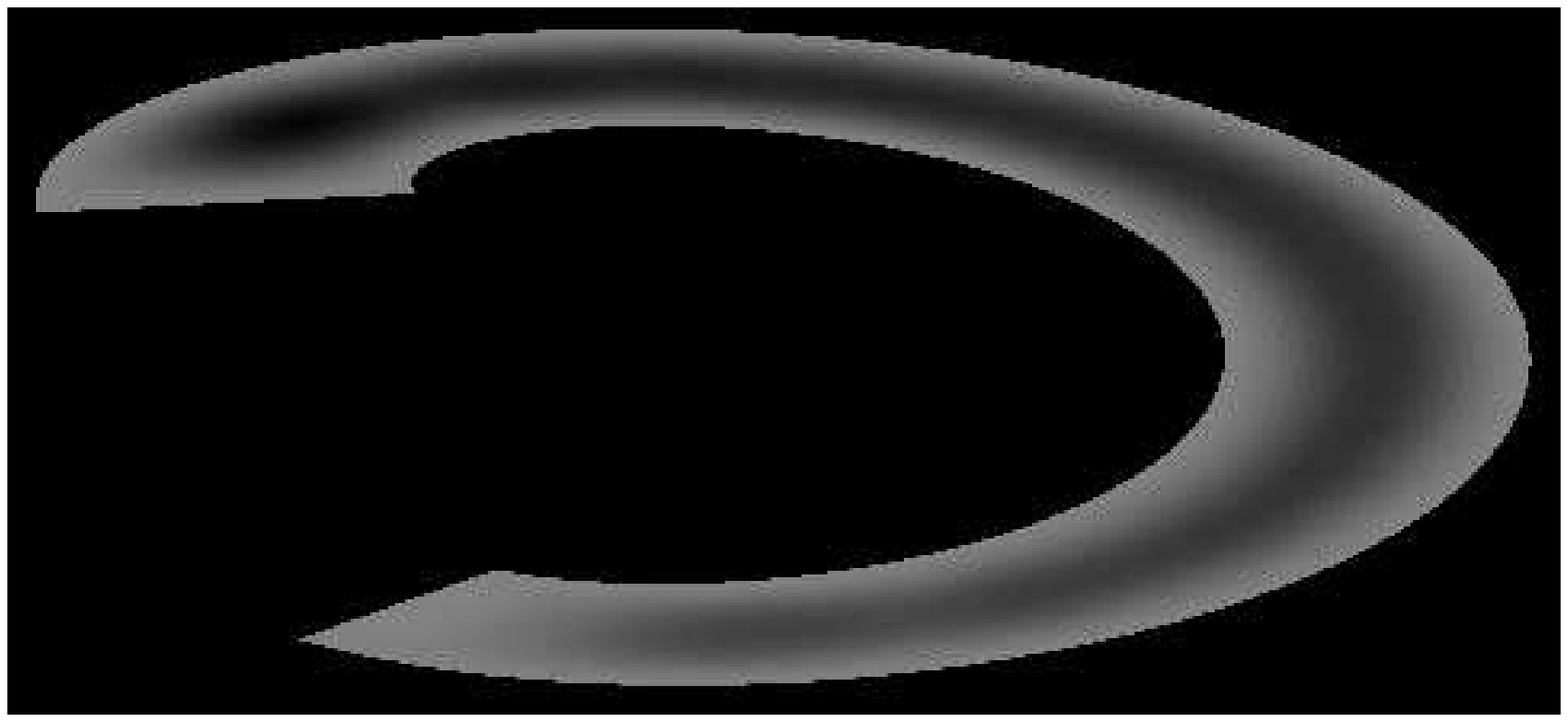} 
&
\epsfxsize=2.5in
\epsffile{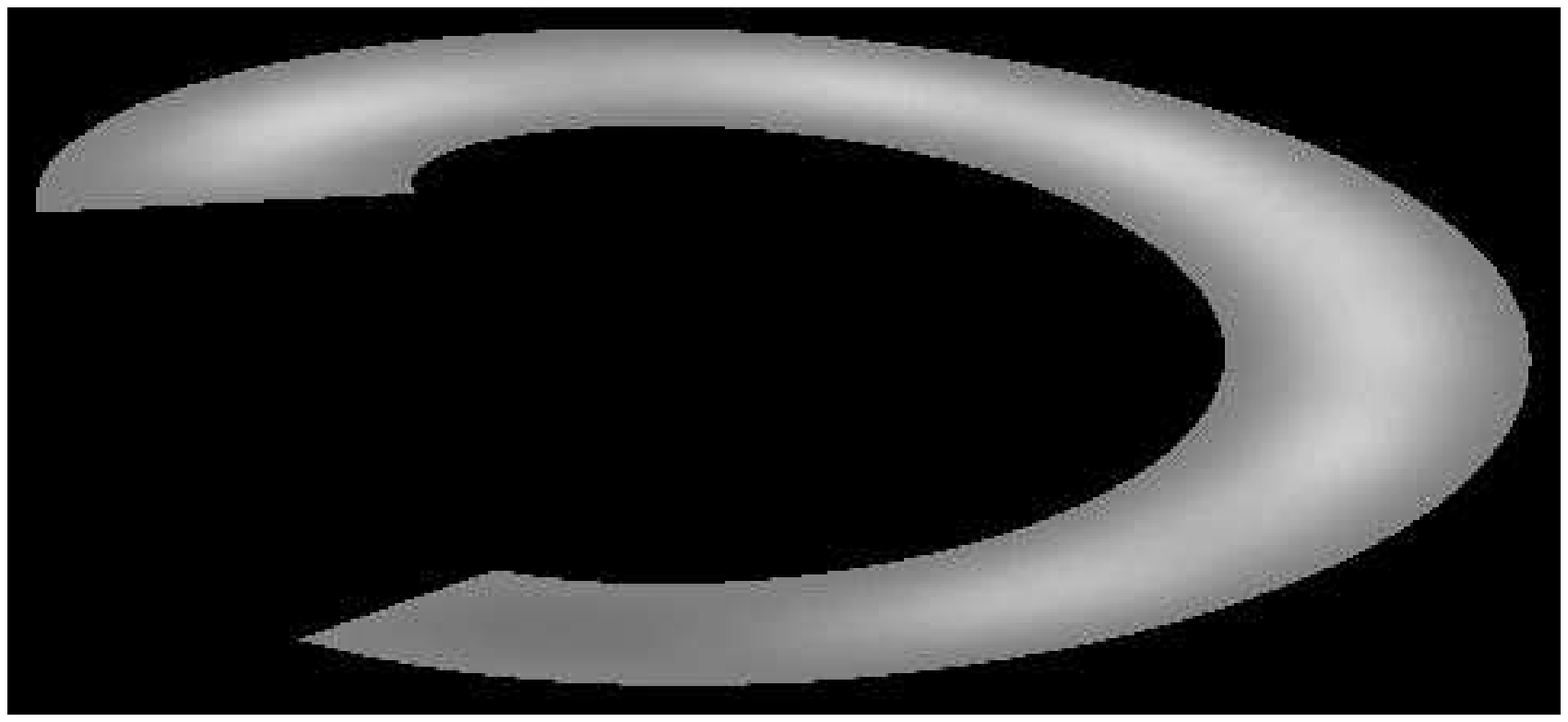} 
\\ % [0.4cm]
\vspace{-3cm}
\epsfxsize=2.5in
\epsffile{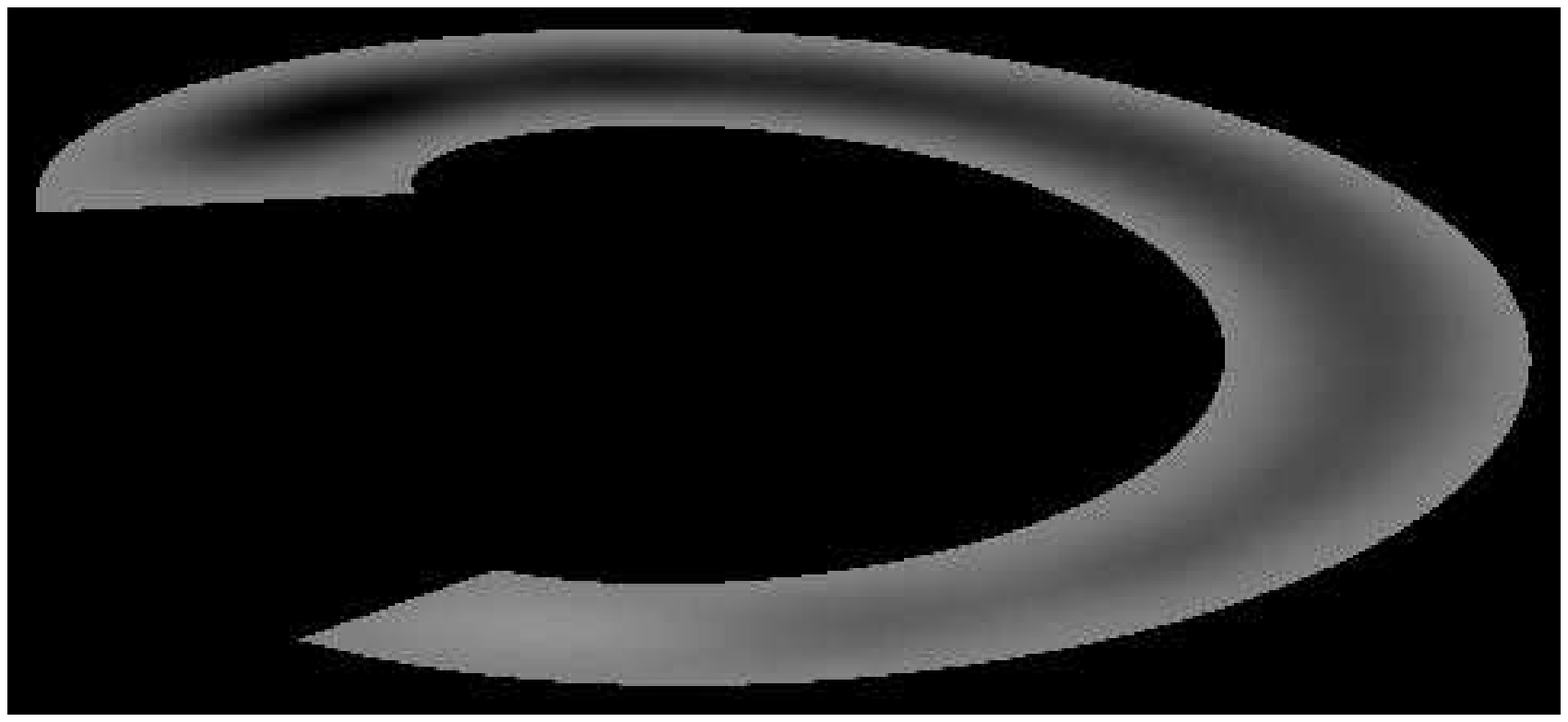} 
&
\epsfxsize=2.5in
\epsffile{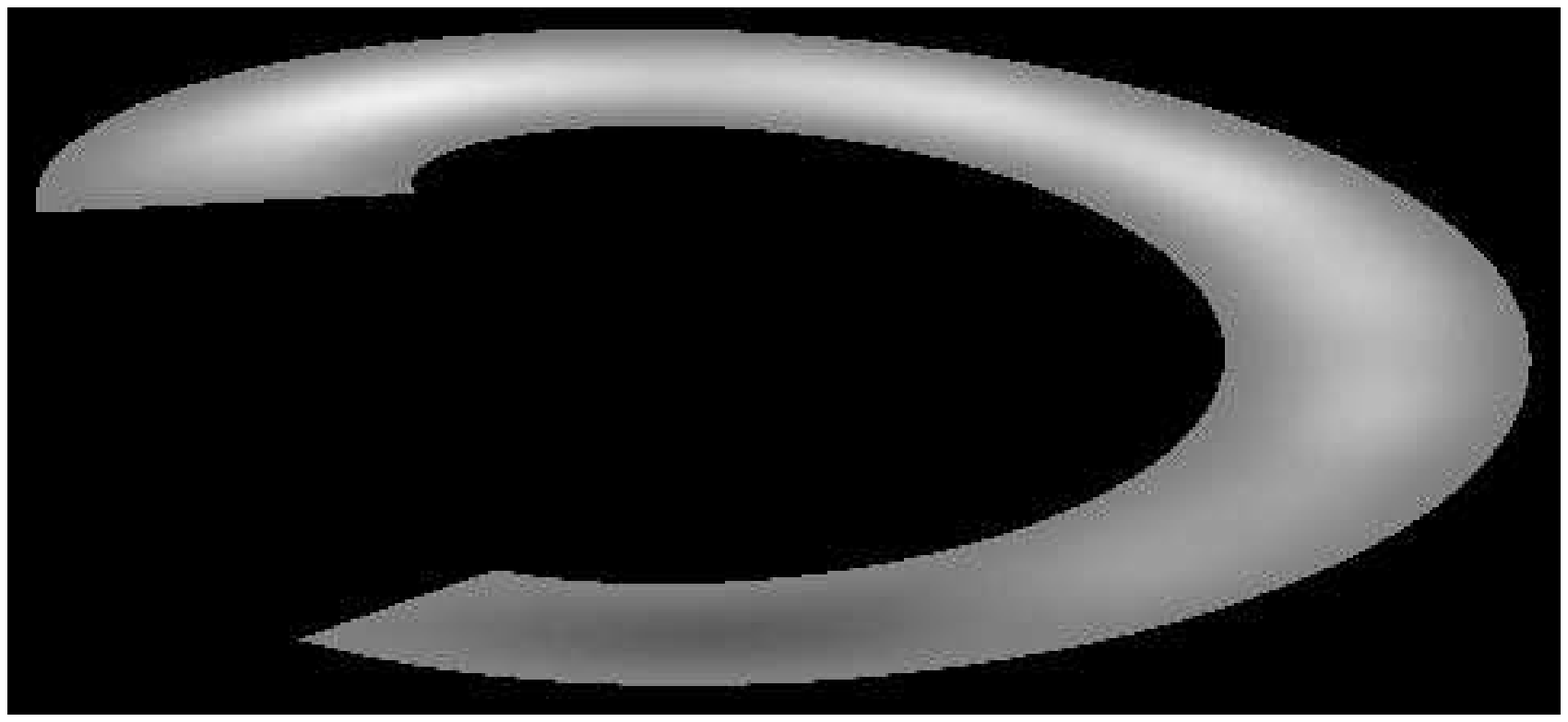} 
\\ % [0.4cm]
\vspace{-3cm}
\epsfxsize=2.5in
\epsffile{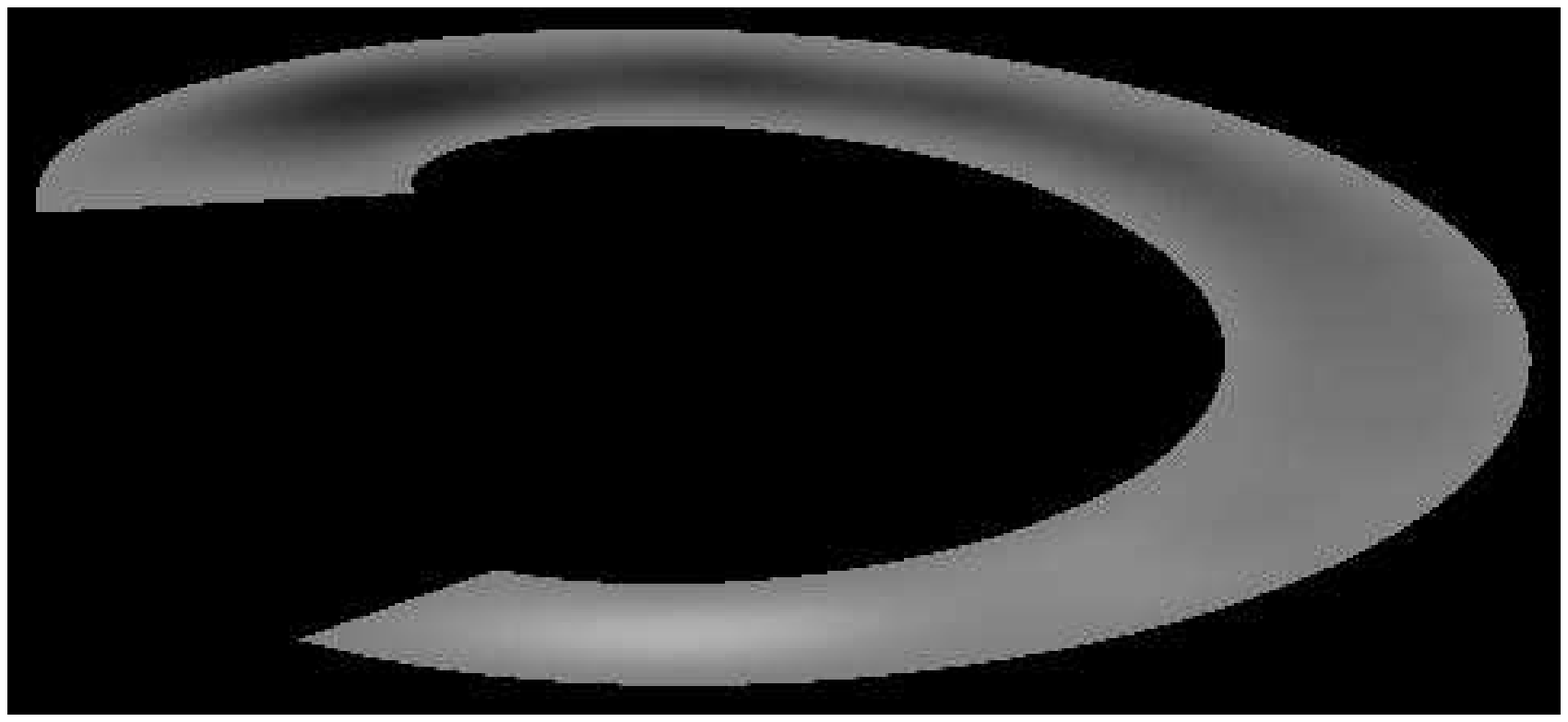} 
&
\epsfxsize=2.5in
\epsffile{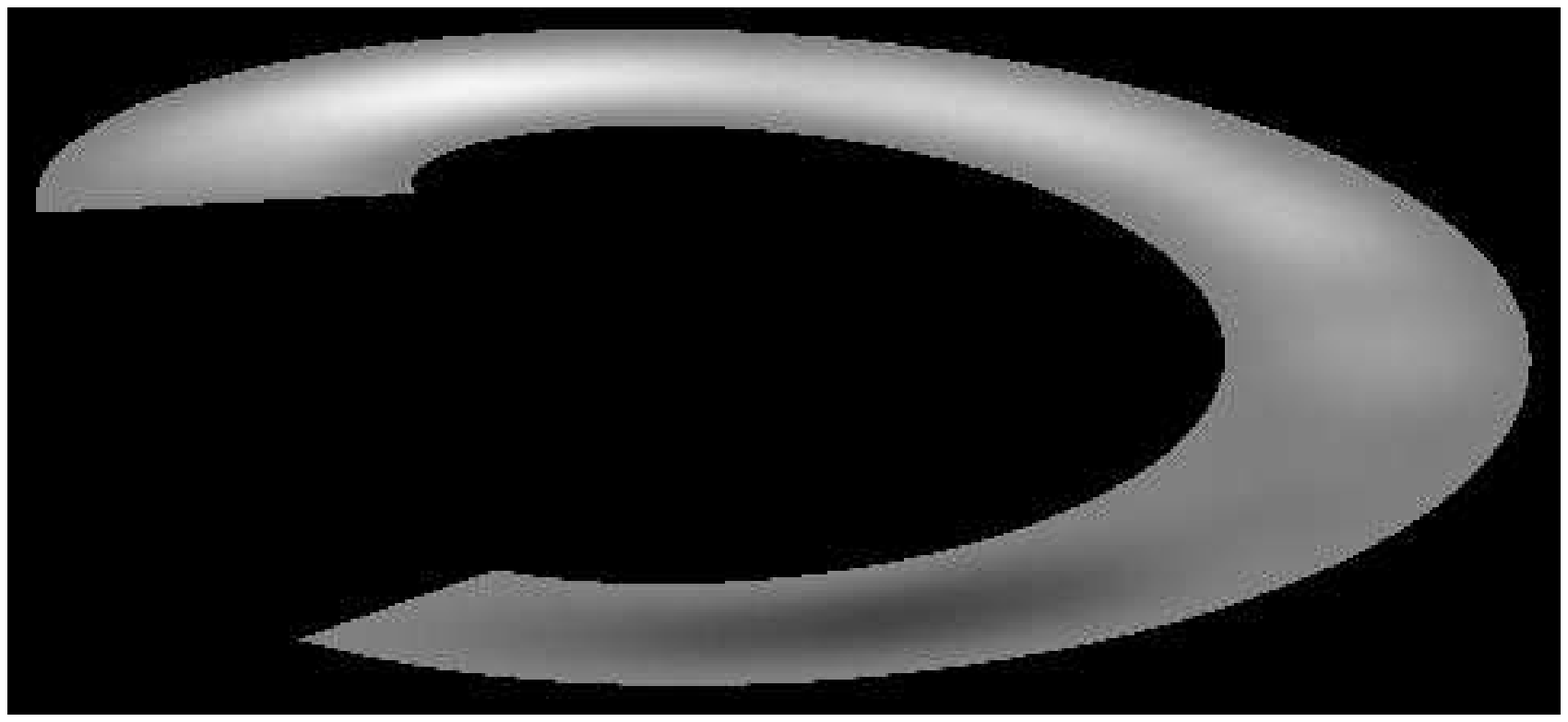} 
\\ % [0.4cm]
\vspace{-3cm}
\epsfxsize=2.5in
\epsffile{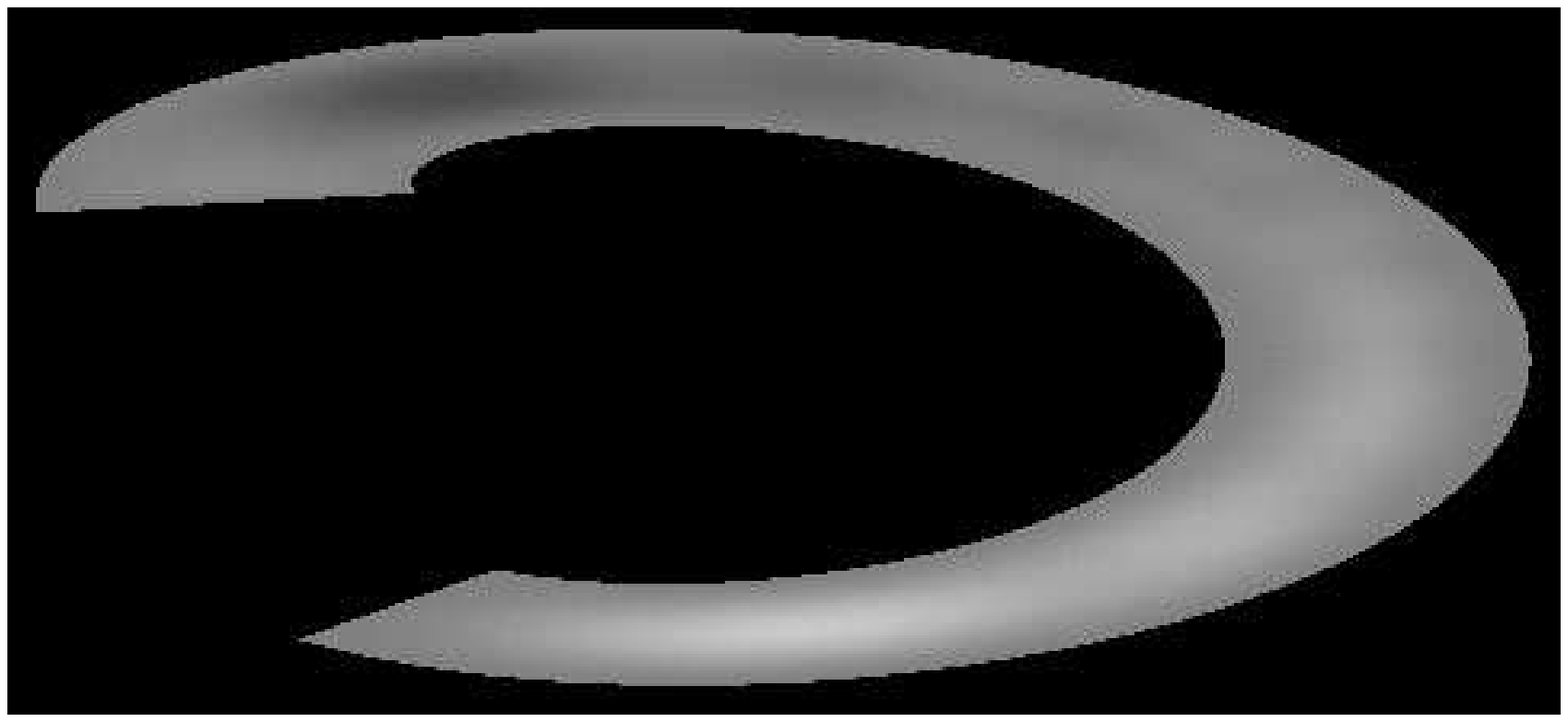} 
&
\epsfxsize=2.5in
\epsffile{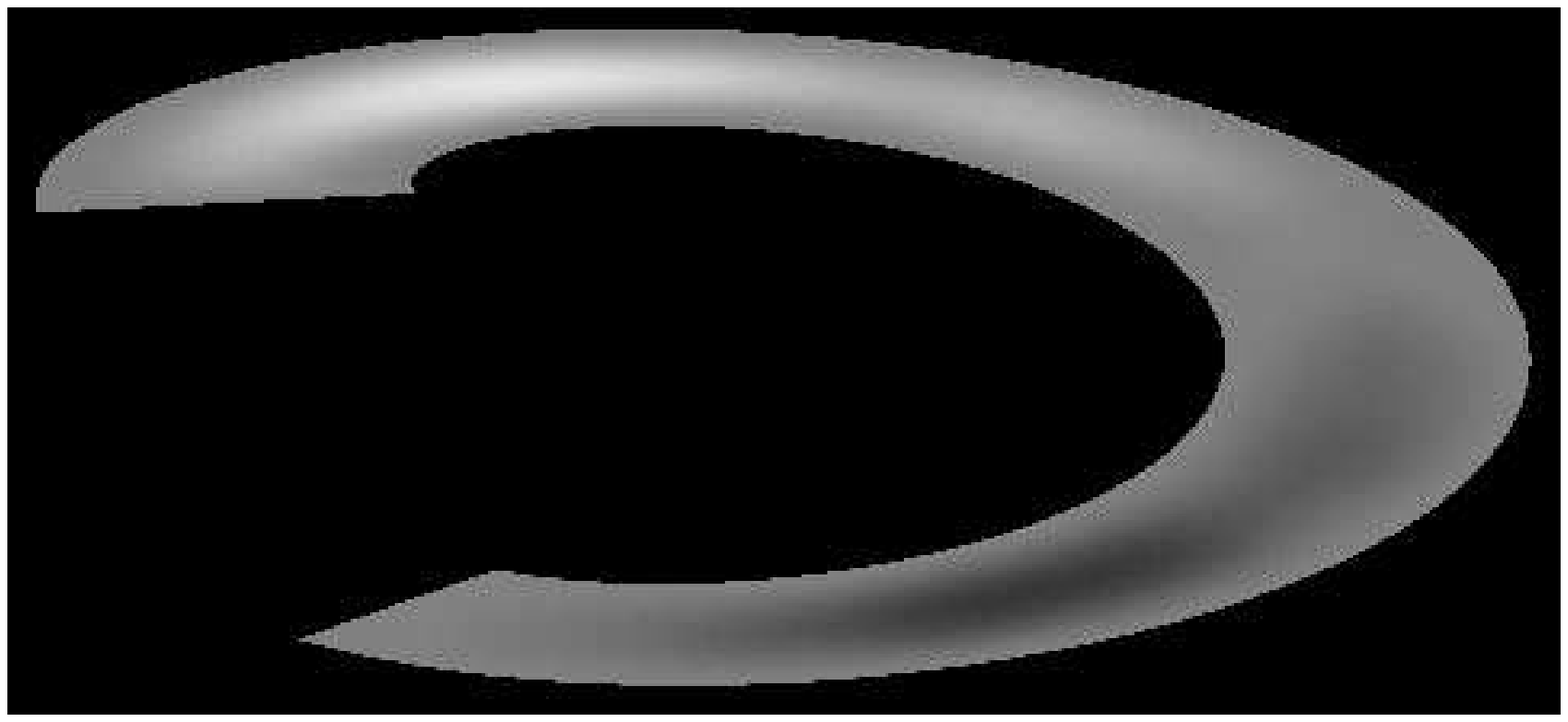} 
\\ % [0.4cm]
\vspace{-3cm}
\epsfxsize=2.5in
\epsffile{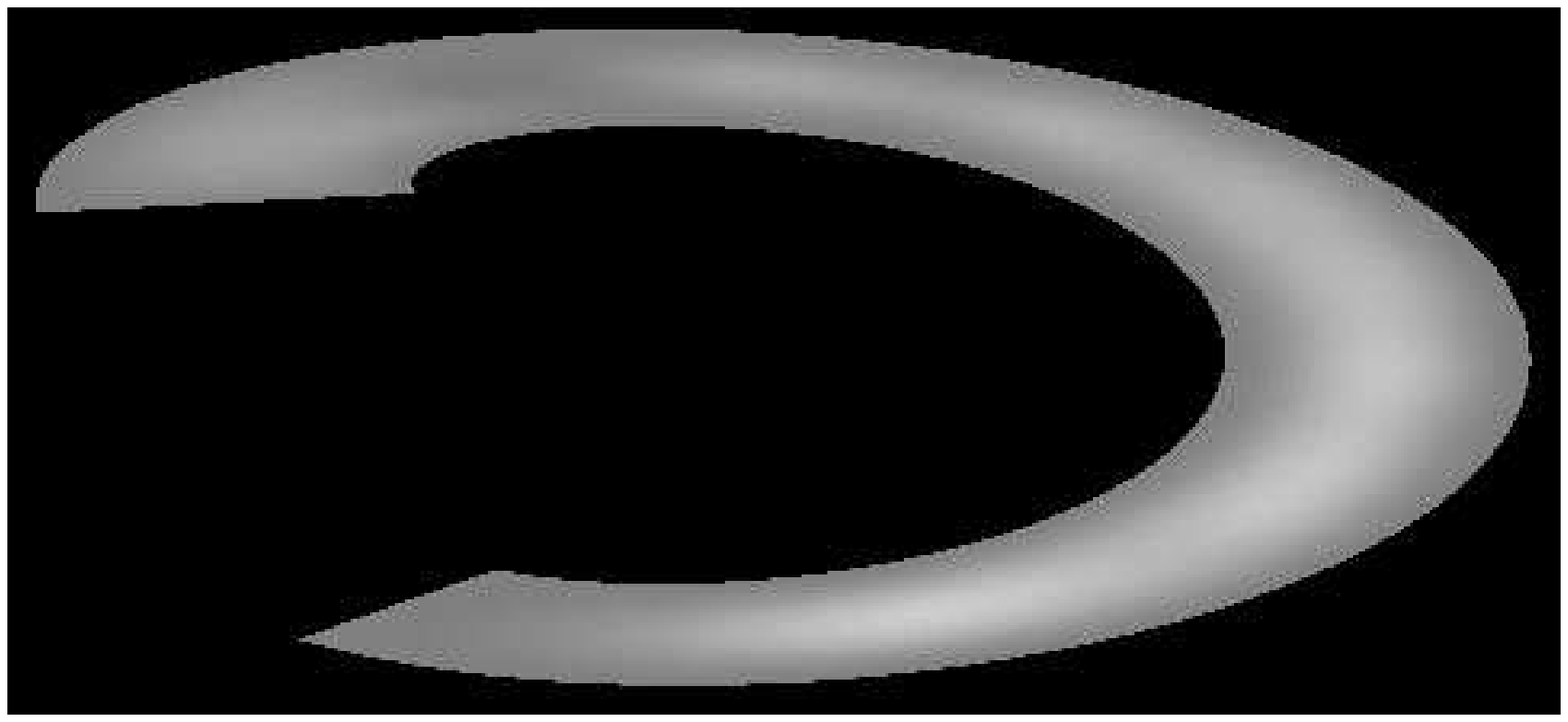} 
&
\epsfxsize=2.5in
\epsffile{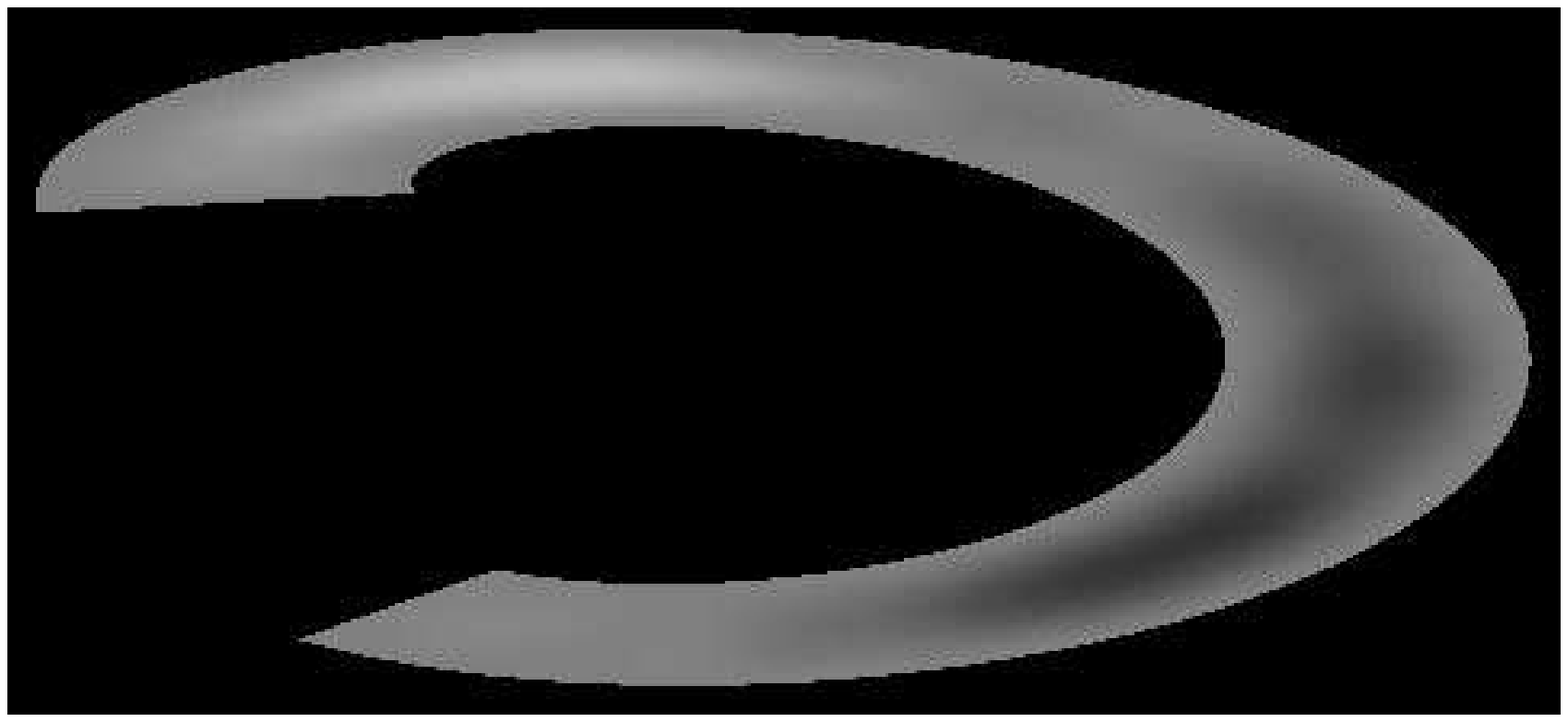} 
\\ % [0.4cm]
% \mbox{\bf (aa)} & \mbox{\bf (bb)}
\end{array}$
\end{center}
\caption{The traveling wave in the test shell.}
\label {fig:wave}
% \caption{The caption for Figure \protect\ref{figtest-fig}}
% \label{figtest-fig}
\end{figure}

In the experiment a traveling wave is produced
% A traveling wave in the shell is observed 
in response to the force impulse in
the fluid.
The wave propagates in the direction of increasing compliance
within the elastic shell.
This agrees with von \Bekesy's observation
that the traveling wave in the cochlea always
propagates from the base to the apex
(see \cite {vonBekesy}).
% As observed by von Bekesy, the wave always
% propagates from the base to the apex
% (see \cite {vonBekesy}).
Snapshots of this wave are  shown in Figure \ref {fig:wave}.
The ten snapshots were taken every 200 time steps beginning 
with time step 600.
The first five appear in the first column, the rest in the second.
Since the displacements of the shell are too small to be seen with the
naked eye, they are represented in the pictures on a gray scale.
Black color indicates the maximal possible displacement down, 
and white, the maximal displacement up (the scale is 
symmetric with respect to zero).
The initial force impulse has pushed the fluid down and,
as can be seen from the first snapshot,
after 600 time steps the shell is displaced downwards.
In the following snapshots we can see that the stiffer part of the shell near
the base is bouncing back and a wave starts propagating towards the apex.
% It is also interesting to note  a distinct two-dimensional pattern
% in the wave.
%
%
%
% \vspace {0.2in}
% \hspace{-0.5in}
%
%================================================================
% \input wave

\section* {Conclusion and Further Research}
\label {chapter:conclusion}

This paper describes an extension of the immersed boundary method
to elastic shells in a viscous incompressible fluid.
The method was developed as a part of a project to construct 
a three-dimensional 
computational model of the cochlea.
It is based on shell equations derived 
using techniques of differential geometry.
The resulting method is a practical method 
which has been implemented and tested on a prototype of a piece of
the basilar membrane.
We have examined the convergence of the algorithm in this case.
Numerical experiments indicate that the algorithm has the first order
convergence rate when the time step is chosen to be linearly
proportional to the fluid mesh width.

The numerical experiments have shown
a traveling wave propagating in the test model shell
in the direction of increasing material compliance
% from the base to the apex of the model shel
in response to external impulsive
excitation of the fluid.
This reproduces the so-called ``travelling wave paradox''
observed in the cochlea: the traveling wave always travels
from the base to the apex of the cochlea 
independently of the location of
the impulse source.

The numerical method described here was subsequently used in
a construction of a complete three-dimensional computational
model of the macro-mechanics of the cochlea.
Additional difficulties in the construction of this
model are caused by the large scale of 
the immersed boundary computations required,
and by the presence in the fluid
of several different elastic materials
in addition to the basilar membrane.
The size of a human cochlea is on the order of
1\,cm $\times$ 1\,cm $\times$ 1\,cm,
while the basilar membrane is about 3.6 cm long and is very narrow
(150---560 \micron).
Therefore, larger fluid grid is necessary in the full cochlea
simulation in order to resolve the 
100 \micron\ scale corresponding to the
width of the basilar membrane.
Immersed boundary computations require large scale computing
resources:
a fluid grid of $256 \times 256 \times 128$ points
needed for the full cochlea model
% a fluid grid of $128^3$ points 
requires a significant amount of
computer memory.
The CFL condition and the stability conditions imposed by
the stiffness forces of the immersed boundaries
force a choice of a very
small time step (approximately 50 nano-seconds)
when the fluid mesh width is small.
The convergence study carried out in this paper indicates that
decreasing the fluid mesh width by a factor of two necessitates
a corresponding decrease in the time step by approximately a factor
of two.
This has indeed been further verified in the construction of the
complete cochlea model.

% This is the major challenge in the construction of the computational
% model for the whole cochlea.
% An application of the immersed boundary method to the complete cochlea
% would require a fluid cube whose volume is approximately 1 cm$^3$,
% which is significantly larger than the fluid cube used in the 
% present model.
% As long as the fluid grid size is not increased this results in
% larger mesh width and smaller precision.
% On the other hand the very small time step required by the CFL
% condition may lead to a significant machine precision error as well.
%
% Despite these obstacles 
% the author has recently succeded in constructing a complete model of
% the cochlea based on the method described in this paper.
The construction of the cochlea model has been achieved in collaboration
with Julian Bunn
on Hewlett-Packard computers
at
% The author has recently completed several successful
% large scale numerical experiments with the full three-dimensional
% model of the cochlea which uses
% the immersed boundary method for elastic
% shells described here.
% The computations were carried out on the HP V2500 computer at
the Center for Advanced Computational Research at Caltech.
The results of this work will be described in future publications.
%================================================================
%
% \input {appendix}
%================================================================
%
%
% \bibliographystyle{plain}
\bibliographystyle{plain}
% \bibliography{../cochlea}
\bibliography{cochlea}
\addcontentsline{toc}{chapter}{Bibliography}
\end{document}